\documentclass[titlepage,final,11pt]{article}
\usepackage{theorem}
\usepackage{enumerate}
\usepackage{array}
\usepackage[reqno]{amsmath}
\usepackage{amssymb}
\usepackage{mathtools}
\usepackage{latexsym}
\usepackage{makeidx}
\usepackage{fancybox}
\input epsf.sty
\usepackage[dvips]{graphicx,color}
\usepackage{showlabels}
\usepackage{subcaption}

\newcommand{\drawing}[4]{
	\begin{figure}[!hbt]
	\begin{center}
	\leavevmode
	\epsfxsize=#2
	\epsfbox{#1}
	\caption{\small #3}
	\label{#4}
	\end{center}
	\end{figure}}

\theoremstyle{change}
\newtheorem{proclaim}{PROCLAIM}[section]
\newtheorem{theorem}[proclaim]{Theorem}
\newtheorem{definition}[proclaim]{Definition}

\newtheorem{proposition}[proclaim]{Proposition}
\newtheorem{corollary}[proclaim]{Corollary}
\newtheorem{example}[proclaim]{Example}

\numberwithin{equation}{section}

\outer\def\proclaim #1. #2\par{\medbreak \noindent{\bf#1.\enspace}{\sl#2}\par
  \ifdim\lastskip<\medskipamount
  \removelastskip\penalty55\medskip\fi}
\def\state #1. { \noindent{\bf#1.\enspace}}
\def\algo #1. { \noindent{\bf#1.\enspace}}

\DeclareMathOperator{\dist}{dist}

\DeclareMathOperator{\dom}{dom}

\DeclareMathOperator{\epi}{epi}

\DeclareMathOperator{\nt}{int}

\newcommand{\comp}{\,{\raise 1pt \hbox{$\scriptstyle\circ$}}\,}

\newcommand{\reals}{\mathbb{R}}
\newcommand{\Reals}{\overline{\mathbb{R}}}
\newcommand{\natnums}{{{\rm l} \kern -.13em {\rm N} }}
\newcommand{\nats}{\mathbb{N}}
\newcommand{\snats}{{I\kern -.29em N}}
\newcommand{\rats}{{Q\kern -.64em \raise 1pt \hbox{$\scriptstyle |$}\;\,}}
\newcommand{\srats}
	{{Q\kern -.56em \raise 1.2pt \hbox{$\scriptscriptstyle /$}\,}}
\newcommand{\ints}{Z\kern -.46em Z}

\newcommand{\pluss}{\hskip1pt \raise1pt\vbox{\hrule width6pt \vskip1pt \hrule
                    width6pt} \kern-4pt{\lower1pt\hbox{\vrule height6pt
		    \kern1pt\vrule height6pt}}\hskip5pt}
\newcommand{\eop}
	{\hfill{$\vcenter{\hrule height1pt \hbox{\vrule width1pt height5pt
   	 \kern5pt \vrule width1pt} \hrule height1pt}$} \medskip}
\newcommand{\half}
	{{\raisebox{1pt}{$\frac{1}{2}$}}}

\newcommand{\setd}{{ d \kern -.15em l}}
\newcommand{\hatsetd}{ d \hat{\kern -.15em l }}

\renewcommand{\epsilon}{\varepsilon}
\renewcommand{\phi}{\varphi}

\newcommand{\lset}{\big\lbrace}
\newcommand{\mset}{{\,\big\vert\,}}
\newcommand{\rset}{\big\rbrace}

\newcommand{\tto}{\;{\lower 1pt \hbox{$\rightarrow$}}\kern -12pt
           \hbox{\raise 2.5pt \hbox{$\rightarrow$}}\;}
\newcommand{\overto}[1]{\,{\raise 0pt\hbox{$\rightarrow$}}\kern -9pt
     \hbox{\lower 3pt \hbox{$\scriptscriptstyle#1$}}\hskip6pt}
\newcommand{\underto}[1]{\,{\lower 1pt\hbox{$\rightarrow$}}\kern -9pt
     \hbox{\raise 4pt \hbox{$\,\scriptscriptstyle#1$}}\hskip7pt}
\newcommand{\bigoverto}[1]{{\raise 0pt\hbox{$\,\longrightarrow$}}\kern -16pt
     \hbox{\lower 3pt \hbox{$\scriptscriptstyle#1$}}\hskip4pt}
\newcommand{\bigunderto}[1]{\,{\lower 1pt\hbox{$\longrightarrow$}}\kern -16pt
     \hbox{\raise 4pt \hbox{$\,\scriptscriptstyle#1$}}\hskip6pt}
\newcommand{\bigbigto}[2]{\,{\raise 0pt\hbox{$\,\longrightarrow$}}\kern -16pt
     \hbox{\lower 3pt \hbox{$\scriptscriptstyle#2$}}\kern -10pt
     \hbox{\raise 4pt \hbox{$\,\scriptscriptstyle#1$}}\hskip7pt}
\newcommand{\downto}{{\raise 1pt \hbox{$\scriptscriptstyle \,\searrow\,$}}}
\newcommand{\upto}{{\raise 1pt \hbox{$\scriptscriptstyle \,\nearrow\,$}}}

\newcommand{\notimply}
	{\quad\hbox{$\Longrightarrow \kern -14pt {/}$}\hskip6pt\quad}

\newcommand{\lto}{\,{\lower 1pt\hbox{$\rightarrow$}}\kern -10pt
     \hbox{\raise 4pt \hbox{$\, \scriptstyle l$}}\hskip7pt}
\newcommand{\eto}{\,{\lower 1pt\hbox{$\rightarrow$}}\kern -10.5pt
     \hbox{\raise 4pt \hbox{$\, \scriptstyle e$}}\hskip7pt}
\newcommand{\hto}{\,{\lower 1pt\hbox{$\rightarrow$}}\kern -11pt
     \hbox{\raise 4pt \hbox{$\, \scriptstyle h$}}\hskip7pt}
\newcommand{\pto}{\,{\lower 1pt\hbox{$\rightarrow$}}\kern -11pt
     \hbox{\raise 4.5pt \hbox{$\, \scriptstyle p$}}\hskip7pt}
\newcommand{\cto}{\,{\lower 1pt\hbox{$\rightarrow$}}\kern -11pt
     \hbox{\raise 4pt \hbox{$\, \scriptstyle c$}}\hskip7pt}
\newcommand{\gto}{\,{\lower 1pt\hbox{$\rightarrow$}}\kern -11pt
     \hbox{\raise 4.5pt \hbox{$\, \scriptstyle g$}}\hskip7pt}
\newcommand{\sto}{\,{\lower 1pt\hbox{$\rightarrow$}}\kern -11pt
     \hbox{\raise 4pt \hbox{$\, \scriptstyle s$}}\hskip7pt}
\newcommand{\awto}{\,{\lower 1pt\hbox{$\rightarrow$}}\kern -15pt
     \hbox{\raise 4pt \hbox{$\, \scriptstyle aw$}}\hskip7pt}
\def\Nto{\,{\raise 1pt\hbox{$\rightarrow$}}\kern -13pt
     \hbox{\lower 3pt \hbox{$\, \scriptstyle N$}}\hskip7pt}
\def\Cto{\,{\raise 1pt\hbox{$\rightarrow$}}\kern -14pt
     \hbox{\lower 3pt \hbox{$\, \scriptstyle C$}}\hskip7pt}
\def\fto{\,{\raise 1pt\hbox{$\rightarrow$}}\kern -14pt
     \hbox{\lower 3pt \hbox{$\, \scriptstyle f$}}\hskip7pt}

\newcommand{\low}[1]{{\lower1pt \hbox{$\scriptstyle #1$}}}
\newcommand{\loww}[1]{{\lower2pt \hbox{$\scriptstyle #1$}}}
\newcommand{\high}[1]{{\raise1pt \hbox{$\scriptstyle #1$}}}

\newcommand{\nsum}{\mathop{\sum}\nolimits}

\newcommand{\nliminf}{\mathop{\rm liminf}\nolimits}

\newcommand{\nlimsup}{\mathop{\rm limsup}\nolimits}
\newcommand{\ninf}{\mathop{\rm inf}\nolimits}
\newcommand{\nsup}{\mathop{\rm sup}\nolimits}
\newcommand{\nmax}{\mathop{\rm max}\nolimits}

\newcommand{\nnmin}{\mathop{\rm minimize}}
\newcommand{\nnmax}{\mathop{\rm maximize}}

\newcommand{\nargmax}{\mathop{\rm argmax}\nolimits}
\newcommand{\nargmin}{\mathop{\rm argmin}\nolimits}

\newcommand{\bfxi}{\mbox{\boldmath $\xi$}}

\newcommand{\lwdy}[2]{\mathrel{\mathop
        {\raisebox{0.1ex}{\null$#1$}}{\hbox{\kern -1.0em
	{\raisebox{-0.8ex}{$\scriptstyle{\;\to #2}$}}}}}}
\newcommand{\lwwdy}[2]{\mathrel{\mathop
        {\raisebox{0.2ex}{\null$#1$}}{\hbox{\kern -1.0em
	{\raisebox{-1.1ex}{$\scriptstyle{\;\to #2}$}}}}}}
\newcommand{\slwdy}[2]{\scriptsize{{\mathrel{\mathop
        {\raisebox{0.1ex}{\null$#1$}}{\hbox{\kern -1.0em
	{\raisebox{-0.8ex}{$\scriptstyle{\;\to #2}$}}}}}}}}
\newcommand{\slwwdy}[2]{\scriptsize{{\mathrel{\mathop
        {\raisebox{0.2ex}{\null$#1$}}{\hbox{\kern -1.0em
	{\raisebox{-1.1ex}{$\scriptstyle{\;\to #2}$}}}}}}}}

\definecolor{lightgray}{gray}{0.75}
\definecolor{myred}{rgb}{0.55,0,0}
\definecolor{myblue}{rgb}{0,0,0.5} 
\definecolor{mygreen}{rgb}{0,0.5,0} 
\definecolor{purple}{rgb}{0.5,0,0.5} 
\definecolor{turq}{rgb}{0,0.805,0.816} 
\definecolor{maroon}{rgb}{0.51,0,0}
\definecolor{MAROON}{rgb}{0.51,0,0}
\definecolor{redor}{rgb}{0.78,0.078,0.078}
\definecolor{dgreen}{rgb}{0,0.3,0}

\newcommand{\Ex}{\mathbb{E}}

\newcommand{\bcdot}{\,{\raise .2ex \hbox{$\centerdot$}}\,}

\begin{document}


\begin{center}
\begin{large}
{\bf Good and Bad Optimization Models: Insights from Rockafellians}
\smallskip
\end{large}
\vglue 0.3truecm
\begin{tabular}{c}
  \begin{large} {\sl Johannes O. Royset 
                                  } \end{large} \\
  Operations Research Department\\
  Naval Postgraduate School, Monterey, California\\
  joroyset@nps.edu\\
\end{tabular}

\vskip 0.1truecm

{\bf Date}:\quad \ \today

\end{center}

\vskip 0.3truecm

\noindent {\bf Abstract}. A basic requirement for a mathematical model is often that its solution (output) shouldn't change much if the model's parameters (input) are perturbed. This is important because the
exact values of parameters may not be known and one would like to avoid being misled by an output obtained using incorrect values. Thus, it's rarely enough to address an application by formulating a model, solving the resulting optimization problem and presenting the solution as the answer. One would need to confirm that the model is suitable, i.e., ``good,'' and this can, at least in part, be achieved by considering a family of optimization problems constructed by perturbing parameters as quantified by a Rockafellian function. The resulting sensitivity analysis uncovers troubling situations with unstable solutions, which we referred to as ``bad'' models, and indicates better model formulations. Embedding an actual problem of interest within a family of problems via Rockafellians is also a primary path to optimality conditions as well as computationally attractive, alternative problems, which under ideal circumstances, and when properly tuned, may even furnish the minimum value of the actual problem. The tuning of these alternative problems turns out to be intimately tied to finding multipliers in optimality conditions and thus emerges as a main component of several optimization algorithms. In fact, the tuning amounts to solving certain dual optimization problems. In this tutorial, we'll discuss the opportunities and insights afforded by Rockafellians.

\vskip 0.2truecm


\baselineskip=15pt

\section{Introduction}\label{sec:intro}

An optimization model for a particular real-world problem isn't unique. With numerous alternatives being available to the modeler, which one is better? There's the usual trade-off between a ``large,'' presumably accurate, but computationally costly model, and a ``small,'' coarse model easily solved. There's also the issue of sensitivity to changes in model parameters, which is a main subject of this tutorial. Will our decision change significantly if we perturb the parameters slightly? This is a major concern because the ``true'' values of the parameters may not be known and one wouldn't like to be misled by a solution obtained using incorrect values. In fact, practitioners of optimization know very well that it isn't enough to obtain a solution and present it as {\em the} answer to a decision maker. One would need to explore the effect of changes to model parameters. Show that a recommended course of action isn't a modeling artifact. Convince the decision maker that the model is ``valid'' and produces reasonable results. A systematic approach addressing these issues is to consider a {\em family of optimization problems} constructed by changing parameters of concern. The resulting sensitivity analysis uncovers troubling situations with unstable solutions and helps us to distinguish between ``good'' and ``bad'' models.

Viewing a model of interest as producing a family of optimization problems also gives rise to algorithmic approaches. It enables us to assess whether the model has computationally attractive properties and thus is good in that sense. Optimality conditions, multiplier vectors and relaxations stem from this perspective as well. Under ideal circumstances and when properly tuned, the relaxations may even furnish the minimum value for the actual model. The tuning of these relaxations is in of itself a process of optimization; it amounts to solving certain {\em dual} problems.

We organize our thinking about families of problems using the fundamental concept of a {\em Rockafellian,} which can be traced back to the middle of the last century with the pioneering work of Rockafellar \cite{Rockafellar.63}, Gale \cite{Gale.67} and others. Under the name ``bifunction,'' the concept was formalized in the seminal text Rockafellar \cite[Chapter 29]{Rockafellar.70}, extended to infinite dimensions in Rockafellar \cite{Rockafellar.74} and beyond the convex case in Rockafellar \cite{Rockafellar.85}. The name ``Rockafellian'' appears in Royset and Wets \cite{primer}, with ``perturbation function'' (Zalinescu \cite{Zalinescu.02}) and ``bivariate function'' (Bauschke and Combettes \cite{BauschkeCombettes.11}) also being found in the literature.

\begin{definition}{\rm (Rockafellian).}
For the problem of minimizing $\phi:\reals^n\to [-\infty,\infty]$, we say that $f:\reals^m\times\reals^n\to [-\infty,\infty]$ is a {\em Rockafellian} with {\em anchor} at $\bar u \in \reals^m$ if
\[
f(\bar u, x) = \phi(x)~~~\forall x\in\reals^n.
\]
\end{definition}

In a practical setting, we think of $\phi$ in the definition as the objective function produced by a particular modeling effort. Since the function is permitted to take the value $\infty$, it accounts for constraints implicitly as discussed further in the next section. Of course, a minimizer of $\phi$ is useful and indicates a possible course of action, but it alone fails to indicate the effect of changing parameter values. An associated Rockafellian $f$  specifies explicitly the dependence on $m$ parameters and defines the {\em family of problems}
\[
\Big\{~ \nnmin_{x\in\reals^n} \,f(u, x),~~~ u\in \reals^m~\Big\},
\]
among which the {\em actual problem} of minimizing $\phi$ emerges as
\[
\nnmin_{x\in\reals^n} \,f(\bar u, x).
\]
In this tutorial, we show that the expanded view involving the whole family of problems reveals concerning modeling assumptions underpinning the actual problem, quantifies the effect of changes in the perturbation vector $u$ away from the anchor at $\bar u$, furnishes optimality conditions for the problem with associated computational possibilities and defines relaxations as well as supporting dual problems.

Our scope covers a vast number of models. The actual problem can be convex or nonconvex, smooth or nonsmooth. It can involve constraints, including restrictions to integer values. It can be the result of a complicated modeling process that accounts for uncertainty in various ways. For example, $\phi$ might be an expectation of the form $\Ex[g(\bfxi,x)]$ and then the perturbation vector $u$ could specify modeling assumptions about the probability distribution of the random vector $\bfxi$. We note, however, that the perturbation vector $u$ rarely represents directly an inherently random quantity such as future product demand or environmental condition, which one can't expect to know at the time of decision making, and rather it represents modeling assumptions about that quantity and other factors.

Sensitivity analysis of linear programs is well-known and the ``economic interpretation'' of dual variables in that context is often the crowning accomplishment of many introductory optimization courses. Besides the restriction to linear settings, the classical treatment suffers from a limited array of perturbations. Typically, one would only deal with changes on the right-hand side of the constraints, but one could imagine other possibilities as well. There's an extensive literature on more general, local stability results for optimization and variational problems examining metric regularity and calmness (Ye and Ye \cite{YeYe.97}, Ioffe and Outrata \cite{IoffeOutrata.08}, Penot \cite{Penot.10}), tilt-stability
(Poliquin and Rockafellar \cite{PoliquinRockafellar.98}, Eberhard and Wenczel \cite{EberhardWenczel.12}, Lewis and Zhang \cite{LewisZhang.13}, Drusvyatskiy and Lewis \cite{DrusvyatskiyLewis.13}, Gfrerer and Mordukhovich \cite{GfrererMordukhovich.15}), full-stability (Mordukhovich et al. \cite{MordukhovichRockafellarSarabi.13}), connections with iterative schemes (Klatte et al. \cite{KlatteKrugerKummer.12}) and specifics of nonlinear programming (Hager \cite{Hager.14}); see also the monographs Aubin and Ekeland \cite{AubinEkeland.84}, Rockafellar and Wets \cite{VaAn}, Bonnans and Shapiro \cite{BonnansShapiro.00} and Mordukhovich \cite{Mordukhovich.18} as well as the surveys Pang \cite{Pang.97} and Aze \cite{Aze.03}.
Global sensitivity analysis emerges via the truncated Hausdorff distance in Attouch and Wets \cite{AttouchWets.91,AttouchWets.93b} and Royset \cite{Royset.18,Royset.20b}.  A review of all these concepts goes beyond this tutorial and we focus instead on an introductory analysis based on epi-convergence and the Rockafellar condition for optimality as developed in Royset and Wets \cite{primer}, which in turn relies on Rockafellar and Wets \cite{VaAn}. Even this more limited scope allows us to address a vast array of applications and gain significant insight about strengths and weaknesses of a model.

The tutorial starts in the next section with basic notation and some motivating examples. Section \ref{sec:epi} presents the epigraphical point of view, which is the foundation for subsequent developments. Section \ref{sec:stability} applies the concepts to Rockafellians to determine whether a model is overly sensitive to parametric changes. Section \ref{sec:normal} reviews the fundamental concepts of normal cones and subgradients, which are subsequently utilized in Section \ref{sec:quant} to quantify local changes in minimum values. Optimality conditions emerging from Rockafellians are the subject of Section \ref{sec:optim}. The chapter ends with an extended discussion of algorithmic possibilities stemming from a Rockafellian including those associated with dual problems.

\section{Notation and Examples}\label{sec:notation}

We often represent an optimization problem by just an objective function $\phi:\reals^n\to \Reals$, where the {\em extended real line}
\[
\Reals = \reals \cup \{-\infty, \infty\} = [-\infty,\infty]
\]
allows us to express constraints implicitly by setting $\phi(x) = \infty$ for infeasible $x\in\reals^n$. In this setting, the {\em domain} of $\phi$, denoted by
\[
\dom \phi = \big\{x\in \reals^n~\big|~\phi(x) < \infty\big\},
\]
specifies the feasible set; decisions producing an objective function value of infinity are considered intolerable. If $\phi(x)\in\reals$ for all $x\in\reals^n$, then $\phi$ is {\em real-valued} and this is also specified by writing $\phi:\reals^n\to \reals$. Generally, the {\em minimum value} and {\em set of minimizers} of $\phi:\reals^n\to \Reals$ become
\[
\inf \phi = \inf\big\{\phi(x)~\big|~x\in\reals^n\big\}, ~~~~~~~~~~~ \nargmin \phi = \big\{x\in \dom \phi~\big|~\phi(x) \leq \inf \phi\big\}.
\]
We note that the pathological case with $\phi(x) = \infty$ for all $x\in\reals^n$ has $\inf \phi = \infty$, which in fact is attained at all $x\in \reals^n$. Still, the definition specifies that $\nargmin \phi = \emptyset$ in that case because $\dom \phi = \emptyset$. This is meaningful because we don't want to consider a ``cost'' of infinity to be optimal. In fact, the problem of minimizing $\phi$ is infeasible. We say that $\phi:\reals^n\to \Reals$ is {\em proper} if $\phi(x) > -\infty$ for all $x\in \reals^n$ and $\phi(x)<\infty$ for some $x\in\reals^n$. Thus, a proper objective function rules out an infeasible problem, but $\nargmin \phi$ could still be empty as the example $\phi(x) = \exp(-x)$ illustrates. Near-minimizers exist more broadly. For $\epsilon\geq 0$, the set of {\em near-minimizers} of $\phi$ is
\[
\epsilon\mbox{-}\nargmin \phi = \big\{x\in \dom \phi ~\big|~ \phi(x) \leq \inf \phi + \epsilon\big\}.
\]
If $\epsilon>0$ and $\inf \phi>-\infty$, then this set is nonempty as illustrated in Figure \ref{fig:nearmin}.

\drawing{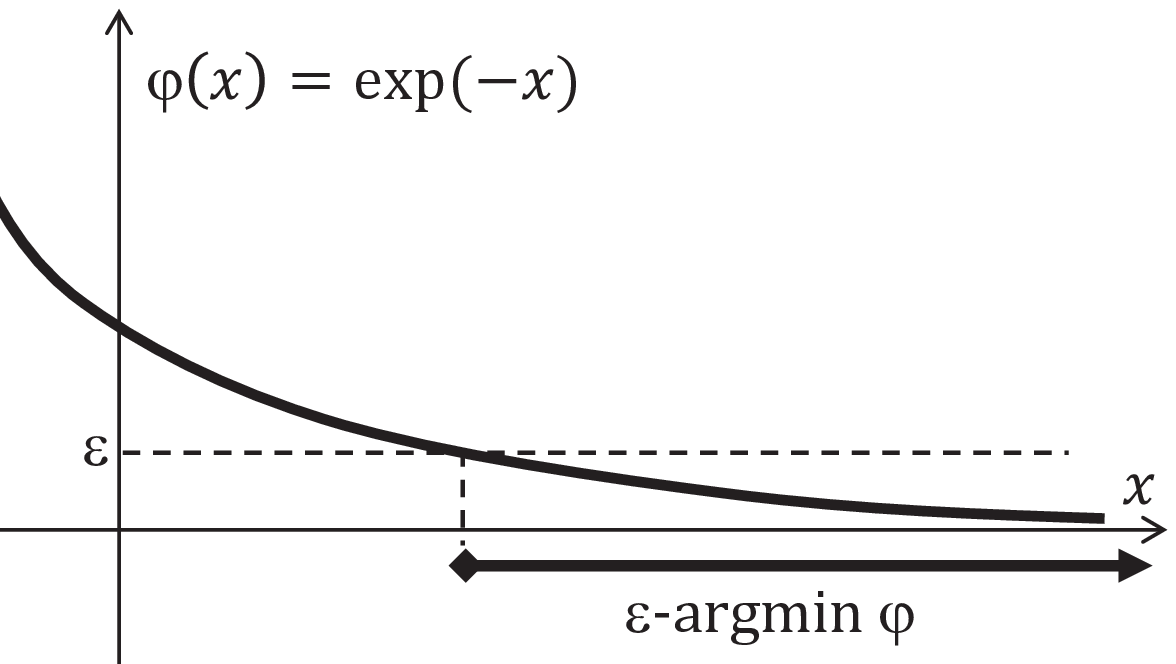}{3.1in}{Example of near-minimizers when $\nargmin \phi$ is empty.}{fig:nearmin}

While it can be useful to hide the details of a model by representing it by a single function, a more detailed analysis and computations should leverage the structural properties of the model. A traditional format is to specify an optimization problem using a real-valued objective function as well as a number of real-valued functions to define equality and inequality constraints. In addition to the glaring imbalance between the detailed specification of the feasible set in terms of {\em many} constraint functions and the lack of details about the objective function, the format fails to bring forth any ``simple'' constraints that should be treated differently than more complicated constraints. From a modeling point of view, there's  also something artificial about singling out one function as objective and the rest as constraints. For example, when designing a system, we might seek a low cost and a low risk of poor performance. To insist on having a single objective function would force us to prioritize cost over risk or vice versa. Moreover, the quantity, being it cost or risk, that has been downgraded to a constraint is subject to a strict requirement of equaling or not exceeding a threshold. In reality, violations might be acceptable if appropriately penalized.

To account for a variety of situations, we often specify optimization models using four components:
\begin{align*}
  & X\subset\reals^n &~~~ \mbox{(basic feasible set)}\\
  & f_0:\reals^n\to \reals &~~~ \mbox{(primary quantity of interest)}\\
  & F:\reals^n\to \reals^m &~~~ \mbox{($m$ secondary quantities of interest)}\\
  & h:\reals^m\to \Reals &~~~ \mbox{(monitoring function for secondary quantities)}
\end{align*}
These components produce an optimization problem in {\em composite form}:
\begin{equation}\label{eqn:compositeform}
    \nnmin_{x\in X} \,f_0(x) + h\big(F(x)\big).
\end{equation}
For example, $f_0(x)$ might be the initial cost associated with decision $x$ and $F(x) = (f_1(x), \dots, f_m(x))$ could be subsequent costs under $m$ different scenarios about future operating conditions each occurring with probability $p_i$. Then, it would be meaningful to set
\[
h(z) = \nsum_{i=1}^m p_i z_i, ~~\mbox{ where } z = (z_1, \dots, z_m),
\]
so that the term $h(F(x))$ becomes the expected subsequent cost across the scenarios.  In another setting, the secondary quantities might be subject to strict requirements expressed by
\[
f_i(x) = 0, ~i=1, \dots, m
\]
or, equivalently, $F(x) = 0$. These requirements are captured by setting
\[
h(z) = \iota_{\{0\}^m}(z),
\]
where $\iota_D$ (Greek letter iota) is the {\em indicator function} of the set $D$ defined as
\[
\iota_D(z) = \begin{cases} 0 & \mbox{ for } z\in D\\
                           \infty & \mbox{ otherwise.}
                           \end{cases}
\]
Thus, $F(x) = 0$ if and only if $\iota_{\{0\}^m}(F(x))<\infty$. Since we always assume that $\alpha + \infty = \infty$ regardless of $\alpha \in \Reals$, a decision $x$ that fails to satisfy the equality requirements would cause $f_0(x) + h(F(x))=\infty$, which indeed implies an infeasible decision. The composite form, especially in the context of certain monitoring functions, is also referred to as ``extended nonlinear programming'' (Rockafellar \cite{Rockafellar.99}).

Let's examine some additional examples in more detail.

\begin{example}{\rm (goal optimization).}\label{eGoalOptimization}
Consider a situation with several quantities of interest (cost, risk, damage, etc.) and the goal of finding a decision that makes all these quantities low. This vague problem statement is a rather common take-away from a first meeting with a prospective client! Of course, there's little hope that one can find a decision that minimizes every quantity at the same time and the goal becomes to identify decisions that balance the various concerns. This leads to the broad area of multi-objective optimization (Miettinen \cite{Miettinen.98}). One approach is then to identify for each quantity of interest a goal $\tau_i$ and a weight $\theta_i$. If $f_i:\reals^n\to \reals$ specifies the values of the $i$th quantity of interest across the decision space, then we can adopt the model
\begin{equation}\label{eqn:goalopt}
\nnmin_{x\in X} \,\sum_{i=1}^m \theta_i \max\big\{0, ~f_i(x) - \tau_i \big\},
\end{equation}
where $X$ is a basic feasible set that restricts the decisions. This is referred to as {\em goal optimization} and falls within the composite form of \eqref{eqn:compositeform}.
\end{example}
\state Detail.  Let $f_0(x) = 0$, $F(x) = (f_1(x), \dots, f_m(x))$ and
\[
h(z) = \nsum_{i=1}^m \theta_i \max\{0,z_i-\tau_i\}
\]
in \eqref{eqn:compositeform} and we recover \eqref{eqn:goalopt}. The resulting model leads to decisions that are as ``close'' as possible to satisfying all the goals. We note that lowering $f_i(x)$ below the goal $\tau_i$ has no benefit, but any value above incurs a per-unit penalty of $\theta_i$.\eop

\drawing{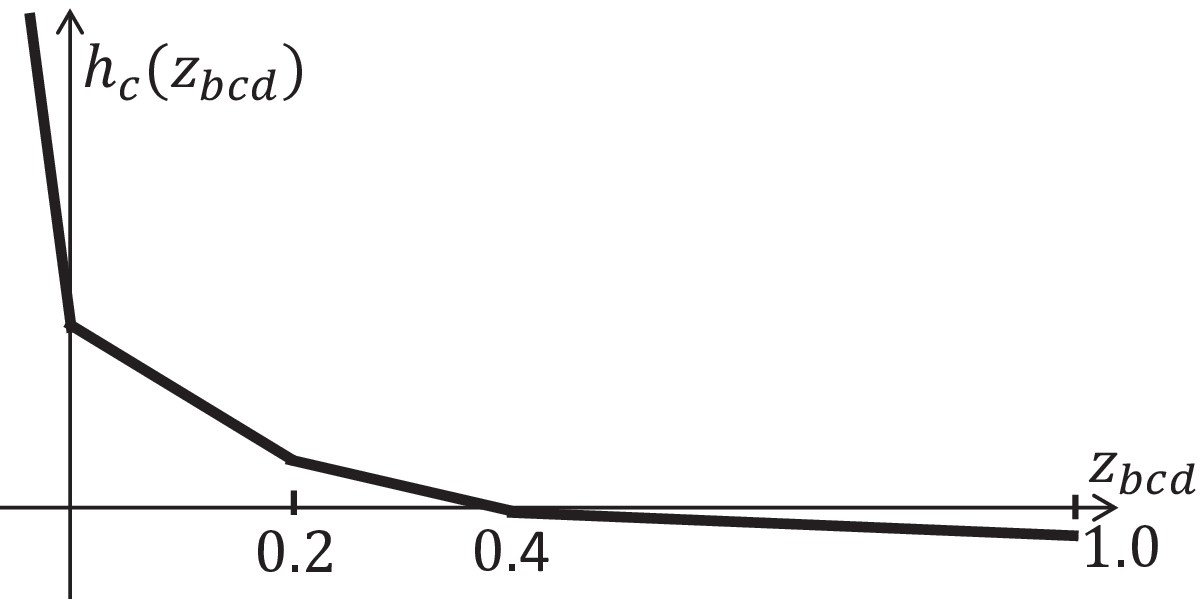}{3.1in}{Monitoring function for inventory level.}{fig:navy}

\begin{example}{\rm (naval resupply at sea).}\label{eNavy}
The US Navy operates worldwide without much need to enter ports for the purpose of resupply. This is achieved by a remarkably low number of transport ships that, as needed, rendezvous with the battle groups at sea. With their numerous moving ``customers,'' the transport ships are faced with a major logistical challenge. Naval planners use an optimization model  to coordinate this effort (Brown and Carlyle \cite{BrownCarlyle.08}). The model tracks the inventory of four commodities (dry stores and food, ship fuel, aircraft fuel and ordnance) for each battle group and one might have expected to find constraints of the kind ``inventory of commodity $c$ on day $d$ for battle group $b$ must be nonnegative.'' However, this isn't the case; transport ships are just encouraged to resupply through a system of penalties and rewards. US Navy ships are allowed to have negative inventory according to the model! This isn't a result of poor modeling. In fact, as we see in this tutorial, there's much merit to the approach. Penalties turn out to produce good models, while (hard) constraints are problematic and easily lead to bad models.
\end{example}
\state Detail. The model in Brown and Carlyle \cite{BrownCarlyle.08} specifies the actions (where to sail, when to rendezvous, what to offload, etc.) for all the transport ships across a planning horizon of up to 180 days. Suppose that the vector $x$ specifies all these actions as well as auxiliary variables and
\[
f_{bcd}(x) \mbox{ is the inventory of commodity $c$ on day $d$ for battle group $b$ under decision $x$,}
\]
which can be determined based on (assumed) known consumption rates. The decision $x$ is constrained by limitations of the transport ships and many other concerns, which we here simply specify by the feasible set $X$. The feasible set, however, doesn't prevent negative values of $f_{bcd}(x)$. The objective function only encourages high values of $f_{bcd}(x)$. Omitting some details, the model takes the composite form \eqref{eqn:compositeform} with $f_0(x) = 0$ and
\[
F(x) = \Big(f_{bcd}(x), ~\forall b,c,d \Big) ~~~~~~~~~~~ h(z) = \nsum_{b,d,c} h_c(z_{bcd}),
\]
where $z$ is the vector with components $z_{bcd}$ and $h_c:\reals\to \reals$ is a monitoring function, which might depend on the commodity $c$; see Figure \ref{fig:navy}. An increasingly large penalty is invoked as the inventory level decreases, while a reward is associated with high levels (above 40\% capacity in Figure \ref{fig:navy}); an inventory above 100\% capacity of the battle group is prevented by constraints in $X$.\eop

Our assessment of optimization models is carried out in the context of certain perturbations as defined by a Rockafellian. For a problem of the composite form \eqref{eqn:compositeform}, one might consider a Rockafellian $f:\reals^m\times\reals^n\to \Reals$ given by
\[
f(u,x) = \iota_X(x) + f_0(x) + h\big(F(x) + u\big),
\]
which then examines the effect of changes to the secondary quantities of interest. Here, the anchor is at $\bar u = 0$ so the actual problem is recovered by minimizing $f(0,x)$ over $x\in\reals^n$. For example, if $h(z) = \iota_{(-\infty,0]^m}(z)$, then the Rockafellian specifies a change from the constraint $F(x) \leq 0$ to $F(x) + u \leq 0$. This corresponds to the typical ``right-hand side'' perturbation of linear programming sensitivity analysis and captures changes to resource budgets. For goal optimization, with $h$ as in Example \ref{eGoalOptimization}, this Rockafellian represents a change in goal for the $i$th quantity of interest from $\tau_i$ to $\tau_i - u_i$.

Any given model can be associated with many different Rockafellians, reflecting concerns about various parameters. The model would be considered good relative to a Rockafellian $f:\reals^m\times\reals^n\to \Reals$ with anchor at $\bar u$ if the Rockafellian exhibits certain desirable properties such as
\begin{quote}
  minimizers and minimum value of $f(u, \cdot\,)$ tend to those of $f(\bar u, \cdot\,)$ as $u\to \bar u$.
\end{quote}
We omit a formal definition and instead give a simple example illustrating that these convergence properties are far from automatically satisfied.

\drawing{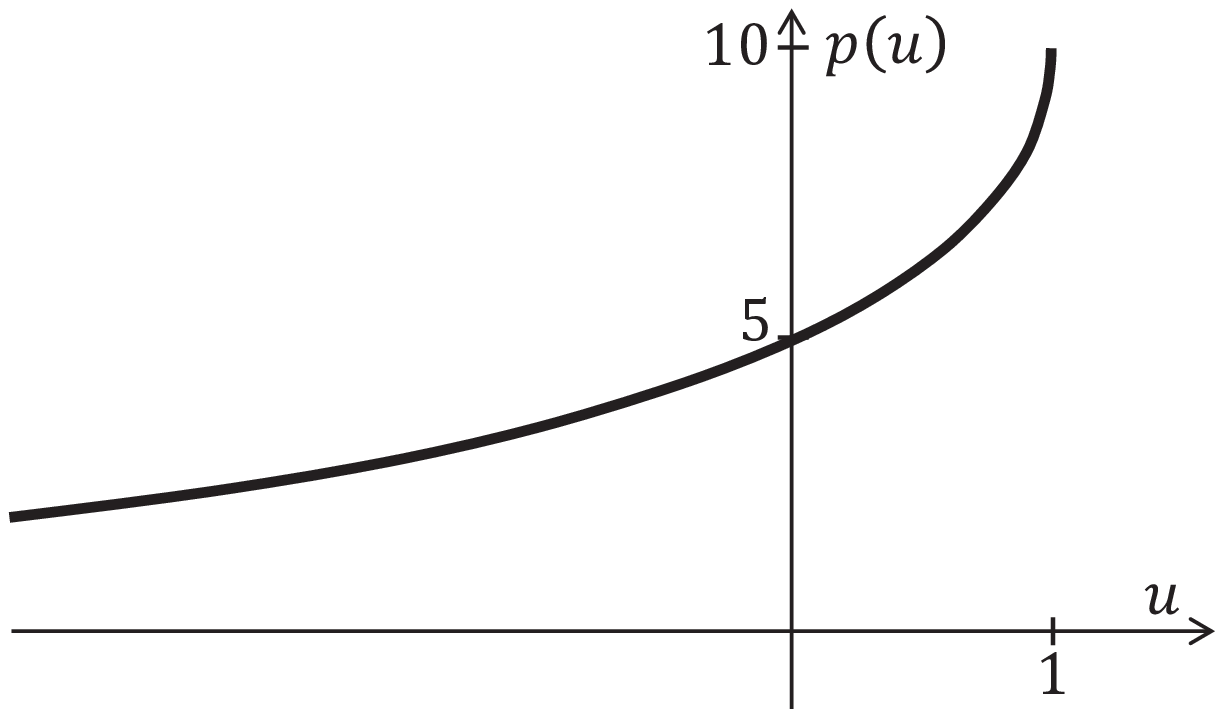}{3.1in}{The minimum value $p(u) = \ninf f(u,\cdot\,)$ as function of $u$ in Example \ref{eSensitivity0}.}{figeSensitivity0}

\begin{example}{\rm (constraint perturbation).}\label{eSensitivity0} The problem of minimizing  $x^2 + 1$ subject to $(x-2)(x-4) + 1$ $\leq 0$ can be associated with a Rockafellian defined by
\[
f(u,x) = x^2 + 1 + \iota_{(-\infty,0]}\big(g(x) + u\big)~~ \mbox{ and }~~ g(x) = (x-2)(x-4),
\]
with anchor at $\bar u = 1$. Let $p(u)=\inf f(u,\cdot\,)$ so that $p(1)$ becomes the minimum value of the actual problem. It turns out that $p(u)$ doesn't tend to $p(1)$ if $u$ approaches $1$ from above. In fact, $p(u) = \infty$ for $u>1$, while $p(1)=10$.  Thus, the Rockafellian highlights the sensitivity to changes on the right-hand side of the constraint in the actual problem. This would have raised concerns about the validity of the model had it represented a real-world situation with constraint parameters being subject to modeling assumptions.
\end{example}
\state Detail.  The feasible set $\dom f(u, \cdot\,)$ is given by the constraint $(x-2)(x-4) + u \leq 0$. A simple application of the quadratic equation produces
\[
\dom f(u, \cdot\,) = \begin{cases}
\big[~3 - \sqrt{1-u}, ~~ 3 + \sqrt{1-u} ~\big] & \mbox{ if } u \leq 1\\
\emptyset & \mbox{ otherwise.}
\end{cases}
\]
As functions of $u$, the minimizers and minimum values become
\begin{align*}
\nargmin f(u,\cdot\,) & = \begin{cases}
\{0\} & ~~~~~\mbox{ if } u < -8\\
\big\{3 - \sqrt{1-u}\big\} & ~~~~~\mbox{ if } -8 \leq u \leq 1\\
\emptyset & ~~~~~\mbox{ otherwise}
\end{cases}\\
p(u) =\inf f(u,\cdot\,) & = \begin{cases}
1 & \mbox{ if } u < -8\\
11 - u - 6\sqrt{1-u} & \mbox{ if } -8 \leq u \leq 1\\
\infty & \mbox{ otherwise.}
\end{cases}
\end{align*}
At every $u<1$, $p$ is continuous. However, at the anchor $\bar u = 1$, the minimum value makes a jump from 10 to $\infty$; see Figure \ref{figeSensitivity0}.\eop

\section{Epigraphs}\label{sec:epi}

The key concept for understanding changes to minimizers and minimum values as well as several other aspects is that of an epigraph. For a function $\phi:\reals^n\to \Reals$, its {\em epigraph} is defined as
\[
\epi \phi = \lset (x,
\alpha)\in \reals^n\times \reals \mset \phi(x) \leq \alpha\rset.
\]
In contrast to the graph of $\phi$, which consists of the points $(x,\alpha)$ with $\phi(x) = \alpha$, the epigraph captures the crucial distinction between $\phi(x) = \infty$ and $\phi(x) = -\infty$ in a manner that's tailored to minimization problems. The former value designates $x$ as infeasible, while the latter value specifies that $x$ must be a minimizer of $\phi$.  Figure \ref{figepi} shows the epigraphs of $f(0,\cdot\,)$, $f(3/4,\cdot\,)$ and $f(1,\cdot\,)$ as defined in Example \ref{eSensitivity0}. These epigraphs appear to converge in the sense that the ``wedge'' $\epi f(u,\cdot\,)$ approaches the vertical line segment $\epi f(1,\cdot\,)$ as $u\to 1$ from below. Moreover, the corresponding minimizers and minimum values converge as already established in Example \ref{eSensitivity0}. This isn't a coincidence. It turns out that convergence of epigraphs essentially guarantees convergence of minimizers and minimum values.

\drawing{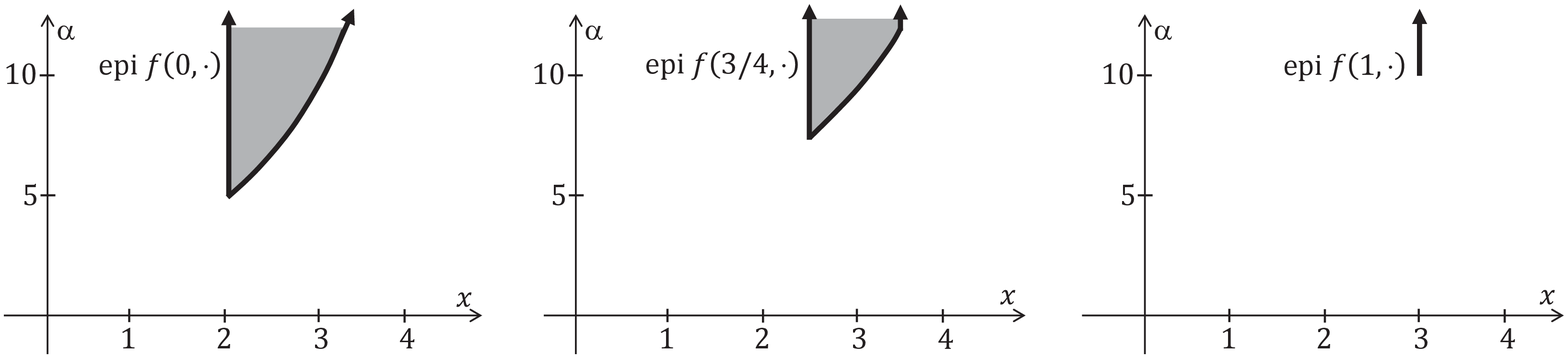}{6.5in}
   {Epigraphs of $f(0,\cdot\,)$, $f(3/4,\cdot\,)$ and $f(1,\cdot\,)$ in Example \ref{eSensitivity0}. Arrows indicate that the epigraphs extend upward indefinitely.}{figepi}

To make precise the meaning of epigraphs converging, we define the {\em point-to-set
distance}\index{distance!point-to-set}\index{dist} between $\bar x\in \reals^n$ and
$C\subset \reals^n$ as the distance between $\bar x$ and the closest point in $C$. Specifically,
\[
  \dist(\bar x,C) = \ninf_{x\in C} \|x-\bar x\|_2 ~~\mbox{ when } ~~C\neq \emptyset~~~~ \mbox{ and } ~~~\dist(\bar x,\emptyset) = \infty.
\]
In the following, we use superscript $\nu$ (Greek nu) to index elements in a sequence, which then runs over the natural numbers $\nats = \{1, 2, \dots\}$.

\begin{definition}{\rm (epi-convergence).}\label{dEpiconv}
For functions $\{\phi, \phi^\nu:\reals^n\to \Reals, ~\nu\in\nats\}$, we say that $\phi^\nu$ {\em epi-converges} to $\phi$, written $\phi^\nu \eto \phi$, when
\[
\epi \phi \mbox{ is a closed set  ~~~and~~ } \dist(z,\epi \phi^\nu) \to \dist(z,\epi \phi)~~\forall z\in \reals^{n+1}.
\]
\end{definition}

Figure \ref{figepi2} illustrates epi-convergence in the context of Example \ref{eSensitivity0}: For any point $z\in \reals^2$,
\[
\dist\big(z, \epi(f(u,\cdot\,)\big) \to \dist \big(z, \epi(f(1,\cdot\,)\big)
\]
when $u\to 1$ from below. The figure visualizes the situation with $z = (1,5)$, in which case one obtains
\[
\dist\Big(z, \epi\big(f(u,\cdot\,)\big)\Big) = \sqrt{ \big(2 - \sqrt{1-u} \big)^2 + \big(6 - u - 6\sqrt{1-u}\big)^2} ~~\mbox{ for } u\in [-8,1].
\]
Certainly, the distance tends to $\sqrt{29}$ as $u\to 1$ from below and this is indeed the distance from $z$ to $\epi f(1,\cdot\,)$. The situation is similar for other $z$.

\drawing{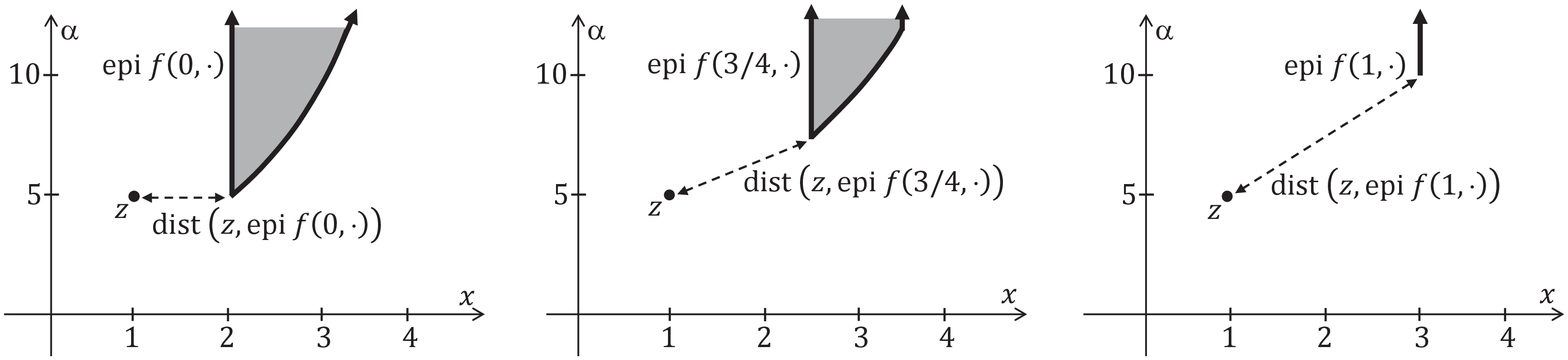}{6.5in}
   {Epi-convergence in Example \ref{eSensitivity0} as indicated by the convergence of the distance from $z$ to the epigraphs.}{figepi2}

The next theorem summarizes the main consequences of epi-convergence; see Theorems 4.14 and 5.5 in Royset and Wets \cite{primer} for a proof and further details.

\begin{theorem}\label{tconvEpiAlgo}{\rm (epigraphical approximations).}
  Suppose that the functions $\phi^\nu:\reals^n\to \Reals$ epi-converge to a proper function $\phi:\reals^n\to \Reals$. Then, the following hold:

  \begin{enumerate}[\rm (a)]
  \item Every cluster point of a sequence $\{x^\nu\in \epsilon^\nu\mbox{-}\nargmin \phi^\nu, ~\nu\in\nats\}$, with $\epsilon^\nu\to 0$, is a minimizer of $\phi$. Thus, in particular, if $x^\nu$ converges to $\bar x$ and each $x^\nu$ is a minimizer of $\phi^\nu$, then $\bar x$ is a minimizer of $\phi$.

  \item Every cluster point $\bar x$ of a sequence $\{x^\nu, \,\nu\in\nats\}$, with $\phi^\nu(x^\nu) \leq \alpha^\nu$ and $\alpha^\nu\to \alpha$, satisfies $\phi(\bar x) \leq \alpha$.

  \item If there's a compact set $B$ such that $B\cap \nargmin \phi^\nu$ is nonempty for all $\nu$, then $\inf \phi^\nu \to \inf \phi$.

  \end{enumerate}
\end{theorem}

We observe that the theorem imposes essentially no assumptions on the functions besides their epi-convergence; convexity and/or smoothness aren't required. Item (b) furnishes the guarantee that a point with a low function value according to the approximation $\phi^\nu$ also has a low value according to the actual function $\phi$.

The following proposition is often helpful in establishing epi-convergence; see Theorem 4.15 in Royset and Wets \cite{primer}. Occasionally, we abbreviate ``{\em there exist(s)}'' by the symbol $\exists$.

\begin{proposition}
\label{tEpiCnvr}{\rm (characterization of epi-convergence)}.
For $\phi,\phi^\nu:\reals^n\to \Reals$, $\phi^\nu\eto \phi$ if and only if the following hold at each $x\in\reals^n$:
\begin{enumerate}[{\rm (a)}]

\item $\forall x^\nu \to x$, one has $\nliminf \phi^\nu(x^\nu) \geq \phi(x)$

\item $\exists x^\nu
    \to x$ such that $\nlimsup \phi^\nu(x^\nu) \leq \phi(x)$.
\end{enumerate}
\end{proposition}

\section{Stability}\label{sec:stability}

Let's now bring in the concept of epi-convergence to examine whether an optimization model is good in the sense that, relative to a particular Rockafellian, minimizers and minimum values change continuously as model parameters are perturbed. The following theorem, which is an immediate consequence of Theorem \ref{tconvEpiAlgo}, shows that the key property to check is epi-convergence of the Rockafellian.

\begin{theorem}{\rm (stability)}.\label{tConsequenceEpiCon} Suppose that $f:\reals^m\times \reals^n\to \Reals$ is a Rockafellian with anchor at $\bar u\in\reals^m$ for the problem
\[
\nnmin_{x\in\reals^n} \,\phi(x).
\]
Let $p(u) = \inf f(u, \cdot\,)$ and $P(u) = \nargmin f(u, \cdot\,)$, i.e., the minimum values and minimizers of $f(u, \cdot\,)$ viewed as functions of $u$. If $f(\bar u, \cdot\,)$ is proper and, for any $u^\nu\to \bar u$, one has $f(u^\nu, \cdot\,) \eto f(\bar u,\cdot\,)$ as well as a compact set $B\subset\reals^n$ such that $B\cap \nargmin f(u^\nu,\cdot\,)$ is nonempty for all $\nu$, then the following hold:
\begin{enumerate}[{\rm (a)}]
  \item $p$ is a continuous function at $\bar u$, i.e., $u^\nu\to \bar u$ implies that $p(u^\nu) \to \inf \phi$

  \item $P$ is outer semicontinuous at $\bar u$, i.e., $u^\nu\to \bar u$ implies that any cluster point $\bar x$ of $\{x^\nu \in P(u^\nu), \,\nu\in\nats\}$ satisfies $\bar x\in \nargmin \phi$.
\end{enumerate}
\end{theorem}

The stability property in the theorem is typically desirable and one might seek to develop models that satisfy the stated assumptions. Let's examine several common cases.

\begin{example}{\rm (regularization as perturbation)}.\label{eRegPerturb}
In data analytics and algorithms relying on proximal point methods, an optimization problem of minimizing a continuous function $f_0:\reals^n\to \reals$ over a nonempty closed set $X$ is often modified by a regularization term involving a penalty parameter. The penalty parameter is rarely fixed. What is the effect of changing it? Specifically, the actual problem is
\[
\nnmin_{x\in X} \,f_0(x) +  \theta_0 r(x),
\]
where $\theta_0\in [0,\infty)$ is the penalty parameter and $r:\reals^n\to [0,\infty)$ is a continuous function; for instance $r(x) = \|x\|_1$ as in lasso regression. With the value of the penalty parameter in focus, we define a Rockafellian $f:\reals\times\reals^n\to \Reals$ for the problem as
\[
f(u,x) = \iota_X(x) + f_0(x) + \theta(u) r(x),
\]
where $\theta(u)$ tends to $\theta(0) = \theta_0$ as $u\to 0$. It turns out that regardless of $\theta_0\in [0,\infty)$, the actual problem is stable in the sense of Theorem \ref{tConsequenceEpiCon} under changes to the penalty parameter as long as a boundedness assumption holds.
\end{example}
\state Detail. Let $u^\nu \to 0$. We use Proposition \ref{tEpiCnvr} to confirm that $f(u^\nu,\cdot\,)\eto f(0,\cdot\,)$. For (a) in the proposition, let $x^\nu\to x$. We consider two cases. If $x\in X$, then $\iota_X(x) = 0$ and
\begin{align*}
\nliminf f(u^\nu,x^\nu) & = \nliminf \big( \iota_X(x^\nu) + f_0(x^\nu) + \theta(u^\nu) r(x^\nu) \big)\\
& \geq \nliminf \big( f_0(x^\nu) + \theta(u^\nu) r(x^\nu) \big) = f_0(x) + \theta(0)r(x) = f(0,x).
\end{align*}
If $x\not\in X$, then $\iota_X(x) = \infty$ and $x^\nu\not\in X$ for all $\nu$ sufficiently large  because $X$ is closed. Thus, $\nliminf f(u^\nu,x^\nu) = \infty$, which matches $f(0,x)$. For (b) in the proposition, let $x$ be arbitrary and take $x^\nu = x$ for all $\nu$. Then,
\[
\nlimsup f(u^\nu,x^\nu) = \nlimsup \big( \iota_X(x) + f_0(x) + \theta(u^\nu) r(x) \big) = \iota_X(x) + f_0(x) + \theta(0) r(x) = f(u,x)
\]
and we conclude that  $f(u^\nu,\cdot\,)\eto f(0,\cdot\,)$. Since $f(\bar u,\cdot\,)$ is proper in view of the nonemptyness of $X$, Theorem \ref{tConsequenceEpiCon} applies as long as one can determine a compact set $B$ that ``reaches'' at least some minimizers in $\nargmin f(u^\nu,\cdot\,)$ regardless of $\nu$. If $X$ is compact, then this is trivially the case but many other possibilities exist as well.\eop

\begin{example}{\rm (risk-averse decision making).}\label{xChangeInRisk} Suppose that we would like to select a decision $x$ from a nonempty compact set $X\subset\reals^n$ such that a cost $g(\xi,x)$ is minimized, where $g:\reals^m\times\reals^n\to \reals$. The difficulty is that the cost depends on $\xi$, which is uncertain and can take values in a finite set $\Xi\subset\reals^m$ each occurring with probability $p_\xi$. One modeling possibility is to select $x$ such that the cost across the worst $(1-\alpha)100\%$ outcomes is minimized, where $\alpha \in (0,1)$ conveys the level of risk averseness. If $\alpha$ is near 1, then $x$ is selected such that the very worst outcomes are made less costly as much as possible. If $\alpha$ is near 0, then the average cost governs the decision. For $\bar\alpha\in (0,1)$, these modeling choices result in the following problem
\[
\nnmin_{x\in X, \gamma\in \reals} \,\gamma + \frac{1}{1-\bar\alpha} \sum_{\xi\in \Xi} p_\xi \max\big\{0, g(\xi,x)-\gamma\big\};
\]
see Section 3.C in Royset and Wets \cite{primer} for further details about such superquantile minimization as well as Rockafellar and Royset \cite{RockafellarRoyset.15b} and Rockafellar and Uryasev \cite{RockafellarUryasev.00} for applications in engineering design and in finance, where it's called CVaR minimization. In practice, it's often difficult to know the right level of risk averseness. Under the assumption that $g(\xi,\cdot\,)$ is continuous for all $\xi\in \Xi$, it turns out that the minimum value for the problem changes continuously as $\alpha$ varies in $(0,1)$. In this sense, the underlying model is good.
\end{example}
\state Detail. Let's express $\alpha = e^u/(1+e^u)$, with $u\in\reals$, and adopt a Rockafellian of the form
\[
f\big(u,(\gamma,x)\big) = \iota_X(x) + \gamma + (1+e^u)\sum_{\xi\in \Xi} p_\xi \max\big\{0, g(\xi,x) - \gamma\big\}.
\]
The actual problem is recovered by minimizing $f(\bar u, \cdot\,)$, where $\bar u = \ln (\bar\alpha/(1-\bar\alpha))$ is the anchor. Let's check the assumptions of the Stability Theorem \ref{tConsequenceEpiCon}. Trivially, $f(\bar u, \cdot\,)$ is proper because $X$ is nonempty. The required epi-convergence can be established using Proposition \ref{tEpiCnvr} as in the previous example and we omit the details.

Let $u^\nu\to \bar u$. The existence of a compact set $B \subset\reals^{1+n}$ such that $B\cap \nargmin f(u^\nu,\cdot\,)$ is nonempty for all $\nu$ is confirmed as follows. For $x^\nu\in X$ and $\gamma^\nu\to \infty$, $f(u^\nu,(\gamma^\nu,x^\nu))\to \infty$ because the term inside the summation is nonnegative. For $x^\nu\in X$ and $\gamma^\nu\to -\infty$, $\min_{\xi\in \Xi, \nu\in\nats} g(\xi,x^\nu)$ is bounded from below because $\Xi$ is a finite set, $X$ is compact and $g(\xi,\cdot\,)$ is continuous. Thus,
\[
\gamma^\nu + (1+e^{u^\nu})\sum_{\xi\in \Xi} p_\xi \max\big\{0, g(\xi,x^\nu) - \gamma^\nu\big\}
\]
involves a first term tending to $-\infty$ and a second term tending to $\infty$, with the second term overpowering the first one due to the coefficient $(1+e^{u^\nu})$, which is greater than one. So, in this case as well, $f(u^\nu,(\gamma^\nu,x^\nu))\to \infty$. We've shown that the distance from the origin to $\nargmin f(u^\nu,\cdot\,)$ can't become arbitrarily large.\eop

\begin{example}{\rm (composite form)}.\label{eComptPert}
Let's consider a problem of the composite form \eqref{eqn:compositeform} with a nonempty closed set $X\subset\reals^n$, continuous $f_0:\reals^n\to \reals$, real-valued $h:\reals^m\to \reals$ and continuous $F:\reals^n\to \reals^m$. We're concerned about changes to $F$ and this leads to a Rockafellian $f:\reals^m\times\reals^n\to \Reals$ defined as
\[
  f(u,x) = \iota_X(x) + f_0(x) + h\big( F(x)+u \big),
\]
with anchor at $0$. Under a boundedness assumption, the Stability Theorem \ref{tConsequenceEpiCon} applies and the underlying model is then good in the sense that minimizers and minimum values vary continuously as specified in the theorem. Since the naval resupply problem in Example \ref{eNavy} satisfies these assumptions, we conclude that the resulting logistical plans are actually rather stable relative to changes in the inventory levels and this can be part of the reason the model has been accepted by US Navy planners.
\end{example}
\state Detail. Since $X$ is nonempty and $h$ is real-valued, $f(0,\cdot\,)$ is proper. Let $u^\nu\to 0$. We use Proposition \ref{tEpiCnvr} to confirm that $f(u^\nu,\cdot\,)\eto f(0,\cdot\,)$. For part (a) of that proposition, let $x^\nu\to x$. We consider two cases. If $x\in X$, then $\iota_X(x) = 0$ and
\begin{align*}
\nliminf f(u^\nu,x^\nu) & = \nliminf \Big( \iota_X(x^\nu) + f_0(x^\nu) + h\big(F(x^\nu) + u^\nu\big)\Big)\\
& \geq \nliminf \Big( f_0(x^\nu) + h\big(F(x^\nu) + u^\nu\big) \Big) = f_0(x) + h\big(F(x)\big) = f(0,x)
\end{align*}
by the continuity of $f_0$, $h$ and $F$. If $x\not\in X$, then $\iota_X(x) = \infty$ and, because $X$ is closed, $x^\nu\not\in X$ as well for all $\nu$ sufficiently large. Thus, $\nliminf f(u^\nu,x^\nu) = \infty$, which matches $f(0,x)$. For (b) in the proposition, let $x$ be arbitrary and take $x^\nu = x$ for all $\nu$. Then,
\[
\nlimsup f(u^\nu,x^\nu) = \nlimsup \Big( \iota_X(x) + f_0(x) + h\big( F(x) + u^\nu\big) \Big) = \iota_X(x) + f_0(x) + h\big(F(x)\big)
\]
again using the continuity of $h$. We've established that $f(u^\nu,\cdot\,)\eto f(0,\cdot\,)$. The existence of a compact set $X$ as required by the Stability Theorem \ref{tConsequenceEpiCon} needs to be verified separately, but this is trivial if $X$ is compact and it holds in other cases as well.

The above analysis relies heavily on the fact that the monitoring function $h$ is real-valued and continuous. While this is acceptable when examining goal optimization (Example \ref{eGoalOptimization}) and naval resupply (Example \ref{eNavy}), it would not hold in the case of (hard) constraints subject to right-hand side perturbations as in Example \ref{eSensitivity0}. In fact, that example shows one can't extend the present results without additional assumptions. This highlights the fundamental advantage of modeling requirements using penalties instead of constraints: under perturbations of the underlying functions, penalties produce good models while constraints might not.\eop

\begin{example}{\rm (penalty approach)}.\label{eConstrPen}
Suppose that the actual problem of interest is
\begin{equation}\label{eqn:constraintpen}
  \nnmin_{x\in \reals^n} \,\phi(x) = \iota_X(x) + f_0(x)  + \iota_{(-\infty,0]^m} \big(F(x)\big),
\end{equation}
where $X$ is nonempty closed and both $f_0:\reals^n\to \reals$ and $F:\reals^n\to \reals^m$ are continuous. As discussed in Examples \ref{eSensitivity0} and \ref{eComptPert}, the underlying model here can be bad in the sense that small changes in $F(x)$ might cause large changes in the minimizers and minimum values. It might be prudent to adjust the modeling of the secondary quantities of interest and adopt the alternative problem
\begin{equation}\label{eqn:constraintpen2}
  \nnmin_{x\in \reals^n}\, \iota_X(x) + f_0(x)  + \theta \sum_{i=1}^m \max\big\{0, f_i(x)\big\},
\end{equation}
where $\theta\in [0,\infty)$ is a penalty parameter and $F(x) = (f_1(x), \dots, f_m(x))$. For any fixed $\theta$, the alternative problem is a special case of that in Example \ref{eComptPert} and therefore comes with the desirable stability properties. However, we can go further and make connections with \eqref{eqn:constraintpen} as well.
\end{example}
\state Detail. Let's consider the Rockafellian given by
\[
f(u,x) = \iota_X(x) + f_0(x)  + \theta(u) \sum_{i=1}^m \max\big\{0, f_i(x) + u_i\big\}
\]
associated with \eqref{eqn:constraintpen2}, where the penalty parameter is also made dependent on $u$. If $u^\nu\to 0$, $\theta(u^\nu)\to \infty$ and
\begin{equation}\label{eqn:penrate}
\theta(u^\nu)\max\big\{0,\nmax_{i=1, \dots, m} u_i^\nu\big\}\to 0~~\mbox{ as } \nu\to \infty,
\end{equation}
then $f(u^\nu,\cdot\,)\eto \phi$. Thus, Proposition \ref{tconvEpiAlgo} can be brought in to verify that minimizers and minimum values of a perturbed version of \eqref{eqn:constraintpen2} converge not only to those of \eqref{eqn:constraintpen2} but also to those of the {\em actual} problem.

To verify the claimed epi-convergence we bring in Proposition \ref{tEpiCnvr}. If $x^\nu\to x$, then $\nliminf f(u^\nu,x^\nu) \geq \phi(x)$ by arguments similar to those in Example \ref{eComptPert}. Let $x^\nu = x$ and consider three cases. (i) If $x\not\in X$, then $\iota_X(x) = \infty$ and $\nlimsup f(u^\nu,x^\nu) = \phi(x)=\infty$. (ii) If $x\in X$ and $f_i(x)>0$ for some $i$, then $\phi(x)=\infty$ and $\nlimsup f(u^\nu,x^\nu)$ can't exceed that value. (iii) If $x\in X$ and $f_i(x) \leq 0$ for all $i$, then
\begin{align*}
\nlimsup f(u^\nu,x^\nu) & = f_0(x) + \sum_{i=1}^m \nlimsup \Big( \theta(u^\nu)  \max\big\{0, f_i(x) + u_i^\nu\big\}\Big)\\
& \leq f_0(x) + \sum_{i=1}^m \nlimsup \Big( \theta(u^\nu) \max\big\{0, u_i^\nu\big\}\Big) = f_0(x) = \phi(x)
\end{align*}
and the claim about epi-convergence follows from Proposition \ref{tEpiCnvr}.

We note that the requirement \eqref{eqn:penrate} has no ramification if $u_i^\nu$ approaches zero from below. However, if $u_i^\nu$ is positive, then one has to make sure that the penalty parameter $\theta(u^\nu)$ grows sufficiently slowly so that the product with $u_i^\nu$ vanishes.\eop

\section{Normal Cones and Subgradients}\label{sec:normal}

While continuity of minimum values as $u$ varies provides indications of a good model, the property fails to quantify how fast the values might change. This raises the question of gradients of minimum value functions, but these might not be defined due to nonsmoothness. Moreover, to develop optimality conditions in a general setting we also need to consider nonsmooth functions; see for example the kinks in Figure \ref{fig:navy}. To address these challenges, we turn to subgradients and supporting normal cones.

\begin{definition}{\rm (normal cone).}\label{dnormalconegeneral} For $\bar x \in C\subset \reals^n$, a vector $v\in\reals^n$ is normal to $C$ at $\bar x$ in the regular sense, or simply a {\em regular normal}, if \footnote{The term $o(\|x-\bar x\|_2)$ for $x\in C$ has the property that $o(\|x-\bar x\|_2)/\|x-\bar x\|_2\to 0$ when $x\in C\to \bar x$ with $x\neq \bar x$.}
\[
\langle v, x-\bar x\rangle \leq o\big(\|x-\bar x\|_2\big) ~~~\mbox{ for } x\in C.
\]
The set of all such regular normal vectors is denoted by $\widehat N_C(\bar x)$.

A vector $v\in\reals^n$ is normal to $C$ at $\bar x$ in the general sense, or simply a {\em normal}, if
\[
v^\nu \to v ~\mbox{ for some }~ v^\nu\in \widehat N_C(x^\nu)~  \mbox{ and }~  x^\nu \in C\to \bar x.
\]
The set of all such normal vectors is $N_C(\bar x)$, the {\em normal cone} of $C$ at $\bar x$.
\end{definition}

Figure \ref{fig:normalcone2} illustrates some common situations. On the left, where the boundary of the set $C$ is smooth at $\bar x$, the regular normal vectors all point in the same direction; recall that $\langle v, x-\bar x\rangle$ is nonpositive when the angle between $v$ and $x-\bar x$ is at least 90 degrees. At the corner point $x'$, the regular normal vectors fan out as it becomes easier to form an angle of at least 90 degrees with points in $C$.

\drawing{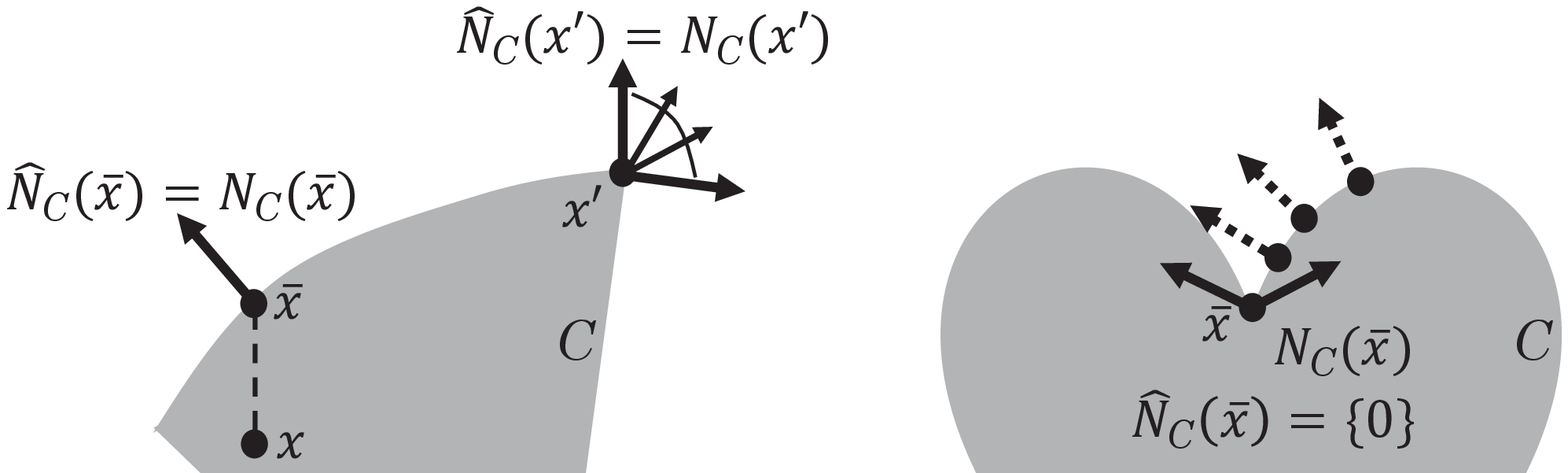}{4.8in}
   {Normal vectors and normal cones.} {fig:normalcone2}

On their own, regular normal vectors fail to provide a solid basis for the treatment of ``complicated'' sets with inward kinks and other irregularities as seen on the right in Figure \ref{fig:normalcone2}. At the inward kink, $\widehat N_C(\bar x)=\{0\}$ since the zero vector is the only vector $v$ that achieves nonpositive $\langle v, x-\bar x\rangle$ for $x$ in $C$ even locally. We achieve a robust notion of ``normality'' by considering the situation at points {\em near} $\bar x$ and this leads to the enrichment of vectors in $N_C(\bar x)$ that aren't regular normals. Figure \ref{fig:normalcone2} shows as dashed arrows regular normal vectors at nearby points approaching $\bar x$ from the right. In the limit, these vectors give rise to the normal vectors at $\bar x$ pointing northwest. Likewise, points approaching $\bar x$ from the left result in the normal vectors at $\bar x$ pointing northeast. In all these illustrations, we think of normal vectors at a point in the sense of ``floating arrows'' and place them relative to that point in $C$. One could just as well visualize normal vectors as forming a cone emanating from the origin.

While normal vectors are important in expressions of optimality conditions, they also define subgradients for arbitrary functions via their epigraphs.

\begin{definition}{\rm (subgradients)}.\label{dSubgrad} For $\phi:\reals^n\to \Reals$ and $\bar x$ with $\phi(\bar x)$ finite, a vector $v$ is a {\em subgradient} of $\phi$ at $\bar x$ if
\[
(v,-1)\in N_{\epi \phi} \big(\bar x, \phi(\bar x)\big).
\]
The set of all subgradients of $\phi$ at $\bar x$ is denoted by $\partial \phi(\bar x)$.
\end{definition}

\drawing{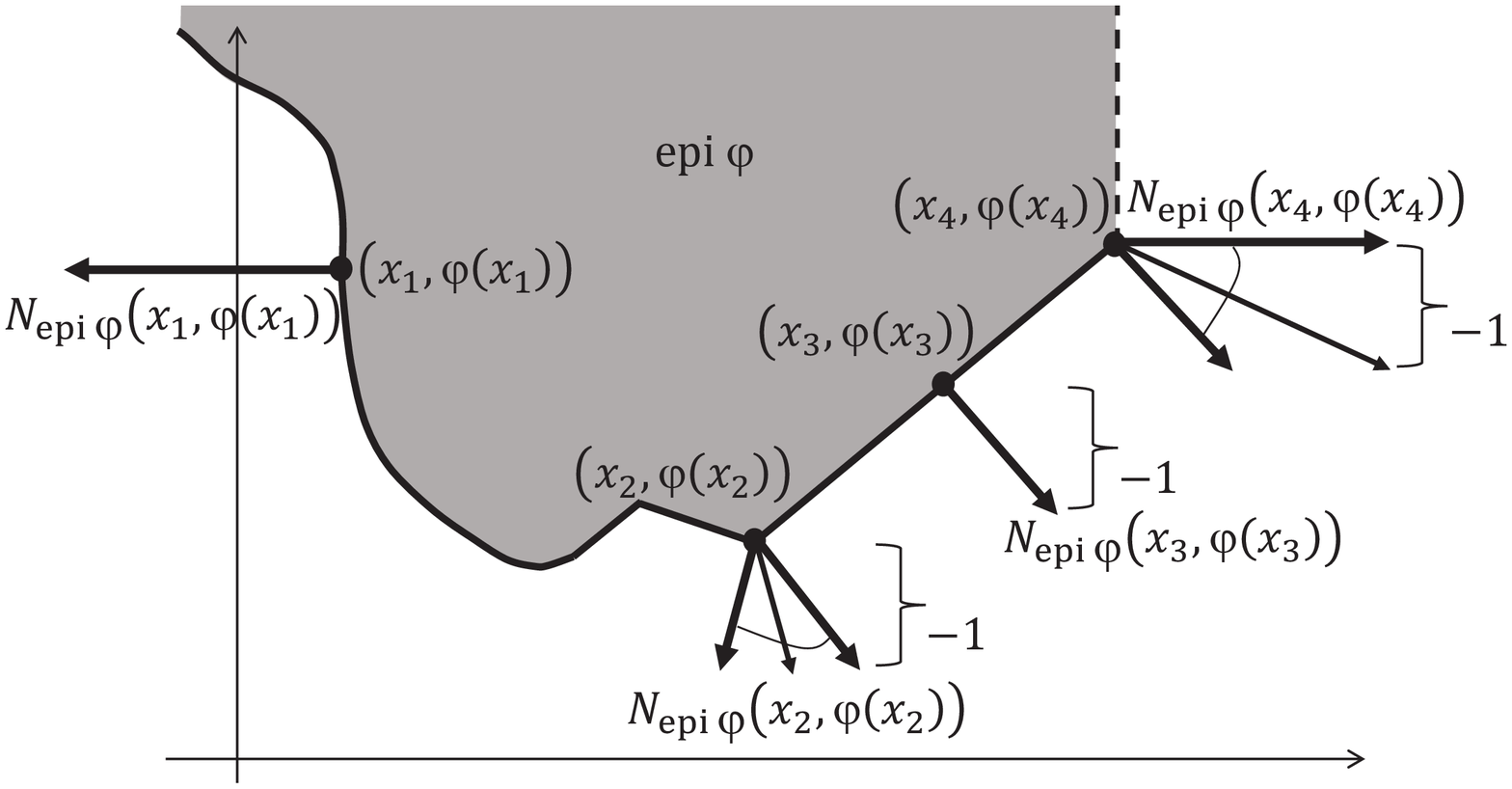}{5.2in} {Normal vectors $(v,-1)$ of epigraphs produce subgradients $v$.} {fig:chainrule0}

Figure \ref{fig:chainrule0} illustrates the normal cones of $\epi \phi$ at four different points. Near $x_3$, the function $\phi$ varies smoothly and the normal vectors at $(x_3, \phi(x_3))$ all point in the same direction, which gives just one vector of the form $(v,-1)$. Thus, there's just one subgradient and it's given by this $v$. The value of this $v$ appears to be roughly 1, which means that $\partial \phi(x_3) = \{1\}$. Of course, this coincides with the slope of $\phi$ at $x_3$. This is confirmed by the first part of the following proposition; under smoothness, subgradients reduce to gradients. A proof is immediate from the definitions.

\begin{proposition}{\rm (subgradients and gradients).}\label{pSubgradGrad2}
For $\phi:\reals^n\to \Reals$ and a point $\bar x$, suppose that $\phi$ is smooth in a neighborhood of $\bar x$. Then,
  \[
  \partial \phi(\bar x) = \big\{\nabla \phi(\bar x)\big\}.
  \]
\end{proposition}

At $(x_2,\phi(x_2))$ in Figure \ref{fig:chainrule0}, the epigraph of $\phi$ has a kink and the normal vectors fan out in many directions. The vectors of the form $(v,-1)$ in the normal cone appear to have $v\in [-1/2, 1]$, which makes this interval the set of subgradients of $\phi$ at $x_2$, i.e., $\partial \phi(x_2) = [-1/2,1]$. The values $-1/2$ and $1$ correspond to the slope of $\phi$ on the left-side and the right-side of $x_2$, respectively.

At points to the right of $x_4$, the value of $\phi$ jumps to infinity. This produces a normal cone of $\epi \phi$ at $(x_4,\phi(x_4))$ that also includes ``horizontal'' vectors, which aren't of the form $(v,-1)$. Still, any normal vector that ``tilts down'' a bit can be made sufficiently long so it becomes in the form $(v,-1)$. Thus, $\partial \phi(x_4) = [1,\infty)$.

At $x_1$, we have the peculiar behavior that $\epi \phi$ is ``vertical'' and all the normal vectors are of the ``horizontal'' kind. Thus, there aren't any normal vectors of the form $(v,-1)$ and $\partial \phi(x_1) = \emptyset$.

As in the case of differential calculus, there's an extensive set of rules that allow us to easily compute subgradients in many practical situations; we refer to Section 4.I of Royset and Wets \cite{primer} for details and simply summarize two key facts here.

\begin{proposition}{\rm (sum rule).}\label{pSubgradGrad3}
For $\phi_1,\phi_2:\reals^n\to \Reals$ and a point $\bar x$, suppose that $\phi_1$ is smooth in a neighborhood of $\bar x$ and $\phi_2$ is finite at $\bar x$. If $\phi$ is the function given by $\phi(x) = \phi_1(x) + \phi_2(x)$, then
  \[
  \partial \phi(\bar x) = \nabla \phi_1(\bar x) + \partial \phi_2(\bar x).
  \]
\end{proposition}

\begin{proposition}{\rm (indicator function).}\label{eSubgradIndicator}
For $X\subset\reals^n$ and $\bar x\in X$, one has
\[
\partial \iota_X(\bar x) = N_X(\bar x).
\]
\end{proposition}

By combining these proposition, we see that the function given by $\phi(x) = f_0(x) + \iota_X(x)$ has
\[
\partial \phi(x) = \nabla f_0(x) + N_X(x)
\]
at any point $x\in X$ provided that $f_0$ is smooth. Regardless of the details of a function, we obtain the following fundamental optimality condition; see Theorem 4.73 of Royset and Wets \cite{primer}.

\begin{theorem}{\rm (Fermat Rule for optimality).}\label{tgenFermat}
For $\phi:\reals^n\to \Reals$ and $x^\star$ with $\phi(x^\star)\in \reals$, one has
\[
x^\star  \mbox{ local minimizer of } \phi ~\Longrightarrow~ 0\in \partial \phi(x^\star).
\]
\end{theorem}

In Figure \ref{fig:chainrule0}, $0\in \partial \phi(x_2) = [-1/2,1]$, which means $x_2$ satisfies the Fermat Rule for optimality. This also illustrates the fact that while a function might be nonsmooth at just ``a few'' points, it tends to occur where we most care: at local minimizers. Consequently, we can't simply ignore the kinks in a function.

The Fermat Rule and our definition of subgradients apply to general functions. Convexity or smoothness isn't needed. In the case of a convex function $\phi:\reals^n\to \Reals$, Definition  \ref{dSubgrad} coincides with the subgradients from convex analysis, which are often stated as those $v\in\reals^n$ with
\begin{equation}\label{eqn:convexinequality}
\phi(x) \geq \phi(\bar x) + \langle v, x-\bar x\rangle ~~~~~\forall x\in \reals^n.
\end{equation}
In the literature, our subgradients are sometimes referred to as ``general,'' ``limiting'' or ``Mordukhovich'' subgradients; see the commentary of Rockafellar and Wets \cite{VaAn} for further details and explanation of why they have emerged as more central than ``Clarke subgradients.''

A special situation occurs when a function becomes ``infinitely steep'' or jumps to infinity at a point as takes place at $x_1$ and $x_4$ in Figure \ref{fig:chainrule0}. The ``horizontal'' normal vectors of the epigraph in such cases are important in dealing with constraint qualifications and related issues, and are given a specific name.

We recall that a function $\phi:\reals^n\to \Reals$ is {\em lower semicontinuous (lsc)} if $\epi \phi$ is a closed set.

\begin{definition}{\rm (horizon subgradients)}.\label{dHorSubgrad} For a lsc function $\phi:\reals^n\to \Reals$ and a point $\bar x$ with $\phi(\bar x)$ finite, a vector $v$ is a {\em horizon subgradient} of $\phi$ at $\bar x$ if
\[
(v,0)\in N_{\epi \phi} \big(\bar x, \phi(\bar x)\big).
\]
The set of all horizon subgradients of $\phi$ at $\bar x$ is denoted by $\partial^\infty \phi(\bar x)$.
\end{definition}

In Figure \ref{fig:chainrule0}, $\partial^\infty f(x_1) = (-\infty, 0]$ and $\partial^\infty f(x_4) = [0,\infty)$. At the points $x_2$ and $x_3$, there aren't any ``horizontal'' normal vectors and thus $\partial^\infty\phi(x_2) = \partial^\infty \phi(x_3) = \{0\}$. The situation at $x_3$ can be settled by the following general fact, which is a direct consequence of the definitions.

\begin{proposition}{\rm (horizon subgradients for smooth functions)}.\label{dHorSubgradSmooth} If  $\phi:\reals^n\to \Reals$ is smooth in a neighborhood of $\bar x$, then $\partial^\infty \phi(\bar x) = \{0\}$.
\end{proposition}

Just as the case for subgradients, horizon subgradients can often be easily computed using a series of calculation rules. We recall two main facts; see Section 4.I of Royset and Wets \cite{primer} for further details.

\begin{proposition}{\rm (expressions for horizon subgradients).}\label{pHorizonSubgradConvex}
If $X\subset\reals^n$ and $\bar x\in X$, then the indicator function $\iota_X$ has
\[
\partial^\infty \iota_X(\bar x) = N_X(\bar x).
\]
If $h:\reals^m\to \Reals$ is proper, lsc and convex and $\bar u\in \dom h$, then
\[
\partial^\infty h(\bar u) = N_{\dom h} (\bar u).
\]
\end{proposition}

A function $\phi:\reals^n\to\Reals$ is {\it locally Lipschitz continuous} at $\bar x$ when there are $\delta \in (0,\infty)$ and $\kappa \in [0,\infty)$ such that
\[
\big|\phi(x) - \phi(x')\big| \leq \kappa\|x-x'\|_2~~~\mbox{ whenever } \|x-\bar x\|_2\leq \delta, ~\|x'-\bar x\|_2\leq \delta.
\]
If $\phi$ is locally Lipschitz continuous at every $\bar x\in\reals^n$, then $\phi$ is locally Lipschitz continuous.

\begin{proposition}{\rm (subgradients and local Lipschitz continuity).}\label{pSubgradLipschitz}
Suppose that  $\phi:\reals^n\to \Reals$ is lsc and $\bar x$ is a point at which $\phi$ is finite. Then,
\[
\phi \mbox{ is locally Lipschitz continuous at } \bar x~~ \Longleftrightarrow ~~ \partial^\infty \phi(\bar x) = \{ 0 \}.
\]
Under these circumstances, $\partial \phi(\bar x)$ is nonempty and compact.
\end{proposition}

\drawing{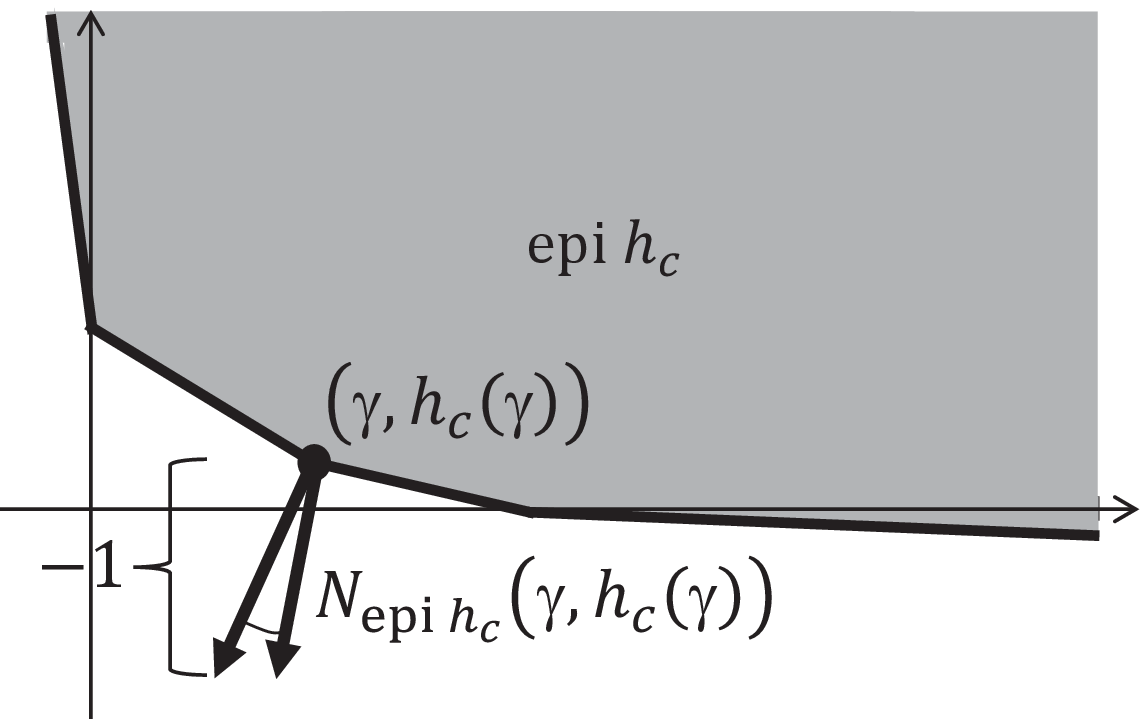}{2.7in} {Normal vectors $(v,-1)$ produce subgradients $v$ of monitoring function in naval resupply problem.} {fig:navyepi}

\begin{example}{\rm (naval resupply at sea; subgradients).}\label{eNavy2} The monitoring function in the
naval resupply problem in Example \ref{eNavy} is expressed by a piecewise affine function of the form
\[
h_c(\gamma) = \max_{i=1, \dots, m} \alpha_i \gamma + \beta_i,
\]
where $\alpha_i, \beta_i$ are given. This function is nonsmooth, but convex. Its subgradients are given by the formula
\[
\partial h_c(\gamma) = \bigg\{ \sum_{i=1}^m \mu_i \alpha_i \bigg| \sum_{i=1}^m \mu_i = 1, ~~\mu_i = 0 \mbox{ if } \alpha_i \gamma + \beta_i < h_c(\gamma), ~~ \mu_i \geq 0 \mbox{ otherwise} \bigg\}.
\]
\end{example}
\state Detail. Figure \ref{fig:navyepi} illustrates the epigraph of $h_c$ and a corresponding normal cone, which in turn produces the subgradients. At most points, the set of subgradients contains only a single value, the slope $\alpha_i$ of the corresponding line segment. However, when the graph of $h_c$ has a kink, the normal cone at that point widens as in the figure and includes all convex combinations of the adjacent slope coefficients.\eop

\section{Quantitative Analysis}\label{sec:quant}

With the confidence that subgradients are defined for arbitrary functions, we now return to the main subject of assessing whether a particular model is good in the sense that the resulting minimum value doesn't change much as a parameter of concern varies. While Section \ref{sec:stability} furnishes conditions under which the minimum value changes continuously, we now step further and quantify the {\em rate of change}. This provides a deeper understanding of the sensitivity to changes in model parameters.  In particular, one might be able to identify which of the parameters are more significant and these can then receive additional scrutiny.

Certain pathological cases might occur if a problem is unbounded in some sense and these are avoided by imposing a boundedness assumption. Specifically, in the context of a function $f:\reals^m\times \reals^n\to \Reals$ with values $f(u,x)$, we say that the function is {\em level-bounded in} $x$ {\em locally uniformly in} $u$ when
\begin{align*}
&\forall \bar u\in\reals^m \mbox{ and } \alpha\in\reals~~~~~\exists \epsilon>0 \mbox{ and a bounded set } B\subset\reals^n\\
&\mbox{such that } \big\{x\in \reals^n~\big|~f(u,x) \leq \alpha\big\} \subset B ~~\mbox{ when } \|u- \bar u\|_2 \leq \epsilon.
\end{align*}
Informally, the property amounts to having, for each $\bar u$ and $\alpha$, a bounded level-set $\{x~|~f(\bar u,x) \leq \alpha\}$ with the bound remaining valid under perturbation around $\bar u$.

Under this assumption, we can, at least partially, characterize the subgradients of the minimum values of a Rockafellian. Since minimum values rarely vary smoothly with parametric changes, we indeed need to bring in  subgradients as gradients may not be defined. The following theorem is given in Theorem 5.13 of Royset and Wets \cite{primer}. For refinements, we refer to Mordukhovich and Nam \cite{MordukhovichNam.05} and Mordukhovich et al. \cite{MordukhovichNamYen.09} and also Aravkin et al. \cite{AravkinBurkeFriedlander.13} for the convex case.

\begin{theorem}{\rm (subgradients of min-value function).}\label{tParametric} For a proper lsc function $f:\reals^m\times\reals^n\to \Reals$, with $f(u,x)$ level-bounded in $x$ locally uniformly in $u$, let
\[
p(u) = \inf f(u,\cdot\,) ~~\mbox{ and }~~ P(u) = \nargmin f(u,\cdot\,) ~~\forall u\in\reals^m.
\]
Then, at any $\bar u\in \dom p$, one has
\[
\partial p(\bar u) \subset \bigcup_{\bar x\in P(\bar u)} \big\{y\in \reals^m~\big|~(y,0)\in \partial f(\bar u, \bar x)\big\}.
\]
If $f$ is convex, then $p$ is convex, the inclusion holds with equality and each set in the union coincides.
\end{theorem}

Typically, $f$ in the theorem is a Rockafellian with anchor at $\bar u$ for an actual problem of minimizing $\phi:\reals^n\to \Reals$. Then, $p(\bar u)$ is the minimum value of the actual problem and $\partial p(\bar u)$ estimates the effect of perturbation on the minimum value. Neither smoothness nor convexity is needed in the first part of the theorem. Though, convexity allows us to sharpen the result from an inclusion to an equality.

In essence, if the subgradients of the Rockafellian are small in length, then the subgradients of the minimum value function are also small. Consequently, the effect on the minimum value of perturbing $u$ away from $\bar u$ is also small. In contrast, if the Rockafellian has large subgradients, then the change in minimum value can be substantial and this calls into question the suitability of the formulation.

\begin{example}{\rm (perturbation of inequality).}\label{eSensitivity} Let's return to Example \ref{eSensitivity0}, but with the constraint in the actual problem now being $g(x)\leq 0$. The Rockafellian then has anchor at $\bar u = 0$. Using Theorem \ref{tParametric}, one can show that
\[
\partial p(0) = \{2\}.
\]
Figure \ref{fig:sensitivity} visualizes the minimum value function $p$ and its affine approximation given by this subgradient, i.e.,
\[
u\mapsto p(0) + 2(u-0) = 5 + 2u.
\]
\end{example}
\state Detail.  We compute the subgradients of the Rockafellian using Proposition \ref{pSubgradGrad3} and this produces
\[
\partial f(u,x) = (0,~2x) + \partial \iota_{(-\infty,0]}\big(g(x) + u\big).
\]
The last term can be calculated using a chain rule (see Theorem 4.64 in Royset and Wets \cite{primer}) so that
\[
\partial f(u,x) = (0, ~2x) + (1, \nabla g(x)) N_{(-\infty, 0]}\big(g(x)+u\big) = (0, ~2x) + (1, ~2x-6)Y(u,x),
\]
for $(u,x)\in \dom f$, where $Y(u,x) = [0, \infty)$ if $g(x)+u = 0$ and $Y(u,x)= \{0\}$ otherwise. Thus,
\[
(y,0)\in \partial f(0,x) ~~ \Longleftrightarrow ~~ 0= 2x+y(2x-6) ~~\mbox{ for some } y\in Y(0,x).
\]
Since we trivially observe that $\bar x = 2$ is the unique minimizer of the actual problem, and this makes $P(\bar u) = \{2\}$ in Theorem \ref{tParametric}, the previous equivalence produces a unique $y = 2$. Consequently,
\[
\partial p(0) \subset \{2\}.
\]
We see directly that $f$ is convex so the inclusion can be replaced by an equality. Moreover, in view of the convexity inequality \eqref{eqn:convexinequality},
\[
p(u) \geq 5 + 2u ~~~~\forall u\in\reals
\]
and this is also confirmed by Figure \ref{fig:sensitivity}.

If we shift the focus to $\bar u = 1$ as in Example \ref{eSensitivity0}, then the above analysis breaks down because there aren't any subgradients of $p$ at $1$; the slope becomes vertical at 1 in Figure \ref{fig:sensitivity}.\eop

\drawing{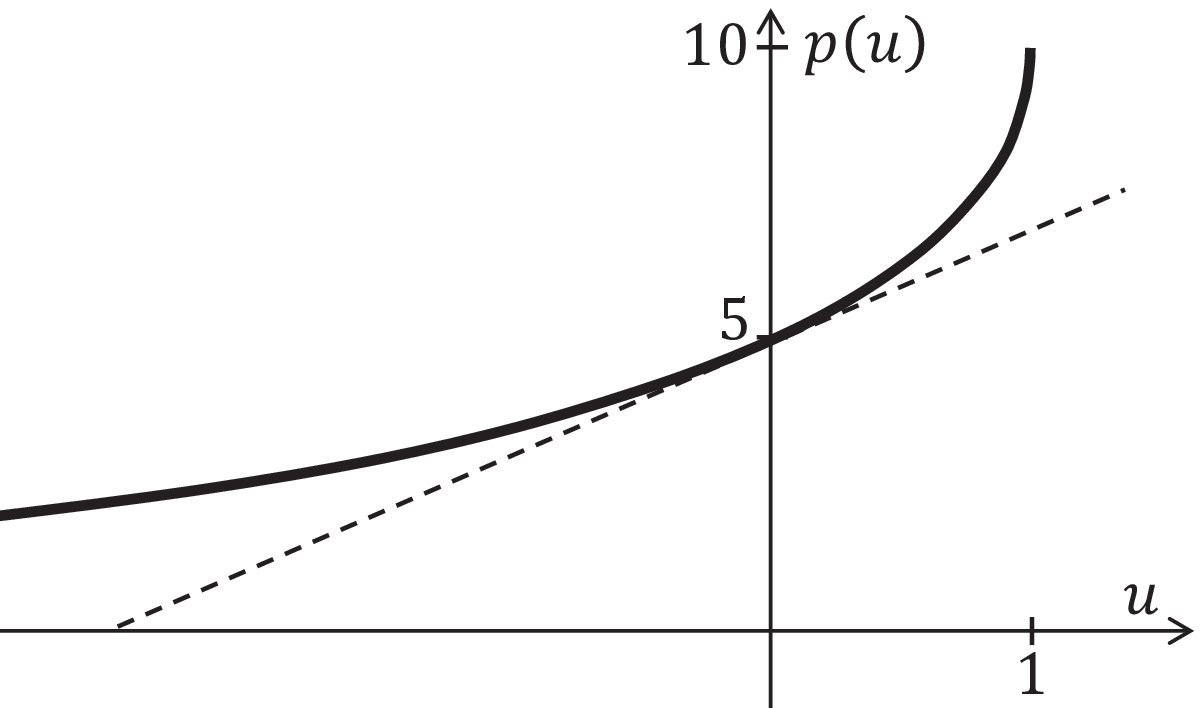}{3.1in} {Minimum value function $p$ and its estimates in Example \ref{eSensitivity}.} {fig:sensitivity}

For problems in the composite form, we obtain the following detailed result; see Proposition 5.16 of Royset and Wets \cite{primer}.

\begin{proposition}{\rm (sensitivity in composite optimization).}\label{pMultAsSubgrad} For smooth $f_0:\reals^n\to \reals$, smooth $F:\reals^n\to \reals^m$, proper, lsc and convex $h:\reals^m\to \Reals$ and nonempty closed $X\subset\reals^n$, consider the problem
\[
\nnmin_{x\in X}\, f_0(x) + h\big(F(x)\big)
\]
and the associated Rockafellian $f:\reals^m\times\reals^n \to \Reals$ given by
\[
f(u,x) = \iota_X(x) + f_0(x) + h\big(F(x)+u\big).
\]
Suppose that $f(u,x)$ is level-bounded in $x$ locally uniformly in $u$. Let $p(u) = \inf f(u,\cdot\,)$, $P(u) = \nargmin f(u,\cdot\,)$ and, for $x\in P(u)$,
\[
Y(u, x) = \big\{y \in \partial h\big(F(x) + u\big)  ~\big|~ -\nabla f_0(x)-\nabla F(x)^\top y\in N_X(x)  \big\}.
\]
Then, for $\bar u\in \dom p$, one has
\[
\partial p(\bar u) \subset \bigcup_{\bar x\in P(\bar u)} Y(\bar u,\bar x).
\]
\end{proposition}

We note that the actual problem in the proposition might be neither smooth nor convex. Still, we obtain a rather specific expression for the subgradients of the minimum value function. In particular, if all the $y$-vectors in $Y(\bar u,\bar x)$ are small, then we conclude that the minimum value changes just a little under perturbation of $u$ away from $\bar u$ and this would indicate a good model formulation.

\begin{example}{\rm (economic order quantity in inventory management).}\index{inventory management} A store manager needs to determine how many units to order of a given product each time the product runs out. The number of units sold per year is $\rho>0$, which means that an order quantity of  $x$ units results in $\rho/x$ orders  per year. Suppose that each one of these orders has a fixed cost of $\beta>0$. (The cost of the units doesn't factor in here because the total number of units ordered across the whole year is always $\rho$.) The manager also faces an inventory cost of $\alpha>0$ per unit and year. With an order quantity $x$, the average inventory is $x/2$ so the annual inventory cost becomes $\alpha x/2$. The manager would like to determine an order quantity $x$ such that both the ordering cost $\beta\rho/x$ and the inventory cost $\alpha x/2$ are low. Let's adopt the goal optimization model (cf. Example \ref{eGoalOptimization})
\[
\nnmin_{x\in X \subset \reals} \,\max\big\{\beta\rho/x - \tau, ~\half\alpha x - \sigma\big\},
\]
where $\tau,\sigma$ are the goals for the ordering cost and the inventory cost, respectively. What's the effect of changing the goals?
\end{example}
\state Detail.  For perturbation $u\in\reals^2$ of the goals, let's adopt the Rockafellian given by
\[
f(u,x) = \iota_X(x)+h\big(F(x) + u\big),
\]
where $h(z) = \max\{z_1, z_2\}$ and
\[
F(x) = \big(f_1(x), f_2(x)\big) = \big(\beta\rho/x - \tau, \half\alpha x - \sigma\big).
\]
Then, $f$ has anchor at $\bar u=0$ and the actual model corresponds to minimizing $f(0,\cdot\,)$. Assuming that $X$ is a nonempty closed subset of the positive numbers, then the assumptions of Proposition \ref{pMultAsSubgrad} are satisfied as long as $f_1$ is extended in a smooth manner from being defined on the positive numbers to all of $\reals$. (This is only a technical issue as we don't consider order quantities below 1 anyway.) We can show that in the notation of Proposition \ref{pMultAsSubgrad},
\begin{align*}
Y(u,x) =  \big\{ &y\in\reals^2~\big|~\beta\rho y_1 x^{-2} - \half\alpha y_2  \in N_X(x)\\
&y_1+y_2 = 1;  ~y_i \geq 0 ~\mbox{ if } f_i(x)+u_i = f(u,x), ~y_i = 0 ~\mbox{ otherwise}; ~i=1,2\big\}
\end{align*}
for $x\in P(u)$. In general, one would employ an algorithm to compute a minimizer $x^\star$ of the model at hand under $\bar u$, which hopefully also computes the corresponding multipliers in $Y(\bar u,x^\star)$. If we assume that $\sigma = \tau=0$ and $X$ isn't active at a minimizer, then we obtain the unique minimizer analytically by solving $f_1(x) = f_2(x)$, which produces
\[
x^\star = \sqrt{2\beta \rho/\alpha}, ~\mbox{ with } Y(0,x^\star) = \big\{\big(\tfrac{1}{2}, \tfrac{1}{2}\big)\big\}.
\]
Then, by Proposition \ref{pMultAsSubgrad}, the minimum value function $p(u)=\inf f(u,\cdot\,)$ actually has $\nabla p(0) =(1/2, 1/2)$. This provides the insight that if the current goal of zero ordering cost is changed to a small positive number $\tau$, which corresponds to setting $u = (-\tau,0)$, then the change in minimum value is approximately $\langle \nabla p(0), u\rangle = -\tau/2$. The negative value is reasonable as raising the goal value reduces the shortfall.
\eop

\section{Optimality Conditions}\label{sec:optim}

The existence of a convenient optimality condition is a major advantage as we analyze and solve an optimization problem. In the nonconvex setting, computational methods might be entirely centered on satisfying the condition. A model can therefore be thought of as having good properties if the resulting optimization problem has an approachable optimality condition. Again, Rockafellians enter as key quantities. The next theorem is a major extension of the Fermat Rule in Theorem \ref{tgenFermat}; see Theorem 5.10 in Royset and Wets \cite{primer}.

\begin{theorem}{\rm (Rockafellar condition for optimality).}\label{tFermatPara}
For the problem of minimizing $\phi:\reals^n\to \Reals$, suppose that $\bar x\in\reals^n$ is a local minimizer, $f:\reals^m\times \reals^n\to \Reals$ is a proper lsc Rockafellian with anchor at $\bar u\in\reals^m$ and the following qualification holds:
\begin{equation}\label{eqn:qualParaFermat}
(y,0) \in \partial^\infty f(\bar u, \bar x) ~~\Longrightarrow~~ y = 0.
\end{equation}
Then,
\[
\exists \bar y\in\reals^m ~\mbox{ such that }~ (\bar y,0) \in \partial f(\bar u,\bar x).
\]
This condition is sufficient for $\bar x$ to be a (global) minimizer of $\phi$ when $f$ is convex.
\end{theorem}

The auxiliary vector $y\in\reals^m$, referred to as a {\em multiplier vector}, is associated with the perturbation vector $u$ and broadens the view of such quantities beyond Lagrange multipliers from the classical Karush-Kuhn-Tucker condition. Every Rockafellian for a problem defines an optimality condition, at least as long as the qualification \eqref{eqn:qualParaFermat} holds. The Rockafellian might now be chosen more with the goal of obtaining a useful optimality condition, than based on concerns about sensitivity analysis. Still, by matching $y$ in Theorem \ref{tParametric} with $y$ here, we realize that solving the Rockafellar condition amounts to determining $(x,y)$, with $y$ furnishing sensitivity information relative to perturbations according to the chosen Rockafellian. This means that solving a problem in this broader sense not only results in a recommended course of action ($x$), but also a measure of stability of the underlying model ($y$) that can help us verify its validity.

For problems in the composite form and the Rockafellian given in Proposition \ref{pMultAsSubgrad}, the optimality condition specializes as stated in Theorem 4.75 of Royset and Wets \cite{primer}.

\begin{theorem}{\rm (optimality for composite problem).}\label{pOptimComposite}
For smooth $f_0:\reals^n\to \reals$, smooth $F:\reals^n\to \reals^m$, closed $X\subset\reals^n$ and proper, lsc and convex $h:\reals^m\to \Reals$, suppose that the following qualification holds at $x^\star$:
\begin{equation}\label{eqn:compositeQual}
y\in N_{\dom h}\big(F(x^\star)\big) ~~ \mbox{ and }~ -\nabla F(x^\star)^\top y\in N_X(x^\star)~~~\Longrightarrow~~~ y=0.
\end{equation}
If $x^\star$ is a local minimizer of the problem
\[
\nnmin_{x\in X} \,f_0(x) + h\big(F(x)\big),
\]
then
\[
\exists y\in \partial h\big(F(x^\star)\big) ~~\mbox{ such that }~ -\nabla f_0(x^\star)-\nabla F(x^\star)^\top y \in N_X(x^\star).
\]
\end{theorem}

The theorem can be specialized further in many directions. We limit the attention to the classical case of equality constraints, which goes back to Lagrange.

\begin{example}\label{eSolProjection}{\rm (smooth objective and equality constraint functions).} For smooth functions $f_i:\reals^n\to \reals$, $i=0, 1, \dots, m$, consider the problem
\[
\nnmin_{x\in\reals^n} \,f_0(x) ~\mbox{ subject to }~ f_i(x) = 0, ~i=1, \dots, m.
\]
If $x^\star$ is a local minimizer with $\{\nabla f_i(x^\star), i=1, \dots, m\}$ linearly independent, then one can apply Theorem \ref{pOptimComposite} to conclude that there is $y\in\reals^m$ such that
\[
\nabla f_0(x^\star) + \nsum_{i=1}^m y_i \nabla f_i(x^\star) = 0, ~~~~f_i(x^\star) = 0, ~~~y_i\in \reals, ~~ i=1, \dots, m.
\]
\end{example}
\state Detail. The problem is of the composite form with $X = \reals^n$, $F(x) = (f_1(x), \dots, f_m(x))$ and $h(z) = \iota_{\{0\}^m}(z)$.  This monitoring function is proper, lsc and convex, with
\[
\partial h(0) = \reals^m
\]
so the multiplier vector $y$ is unrestricted. Since $N_X(x) = \{0\}$, the gradient condition in Theorem \ref{pOptimComposite} reduces to
\[
\nabla f_0(x^\star) + \nabla F(x^\star)^\top y = 0.
\]
The domain of $h$ is $\{0\}^m$. Consequently, $N_{\dom h}(z) = \reals^m$ if $z = 0$ and $N_{\dom h}(z)$ is otherwise not defined. The qualification \eqref{eqn:compositeQual} therefore specializes to checking whether $y\in\reals^m$ and $\nabla F(x^\star)^\top y = 0$ imply $y = 0$. But, this is just the linear independence assumption.\eop

\section{Algorithmic Approaches}\label{sec:algo}

A Rockafellian associated with a particular problem also specifies algorithmic possibilities via problem relaxations. The strength of the resulting relaxations can be used to assess whether the underlying model is computationally attractive in the first place. Moreover, these relaxations define a dual problem, which highlights the many opportunities that lie beyond the classical dual problems of linear and convex programming.

This section starts with the construction of relaxations and then defines the corresponding dual problems.

\subsection{Rockafellian Relaxation}

A {\em relaxation} of a minimization problem is an alternative problem with a minimum value no higher than that of the actual one. To make this concrete, let's consider $\phi:\reals^n\to \Reals$, the actual  problem
\[
\nnmin_{x\in\reals^n} \,\phi(x)
\]
and an associated Rockafellian $f:\reals^m\times\reals^n\to \Reals$ with anchor at $0$. (The focus on $0$ instead of a more general $\bar u$ as the anchor promotes symmetry below, without much loss of generality because one can always shift the perturbation vector by redefining $f$.)  We can't do worse if permitted to {\em optimize} the perturbation vector $u$ together with $x$. Thus, the problem
\[
\nnmin_{u\in\reals^m,x\in\reals^n} \,f(u,x)
\]
is a relaxation of the actual problem. Regardless of $y\in\reals^m$, this is also the case for
\begin{equation}\label{eqn:relax}
\nnmin_{u\in\reals^m,x\in\reals^n} \,f(u,x) - \langle y, u\rangle
\end{equation}
because $u=0$ remains a possibility. In summary, for any $y\in\reals^m$,
\begin{equation}\label{eqn:lagrangianLB}
\ninf_{x\in\reals^n} \phi(x) = \ninf_{x\in\reals^n} f(0,x) \geq \ninf_{(u,x)\in\reals^m\times\reals^n} \big\{f(u,x) - \langle y,u \rangle  \big\}.
\end{equation}
A benefit from bringing in a vector $y$ is that we now can tune the relaxation by selecting $y$ appropriately. We return to this subject below.

We immediately have an approach to addressing the actual problem: (i) Solve the relaxation \eqref{eqn:relax} for some $y\in\reals^m$ and obtain a lower bound on $\inf \phi$. (ii) Use any (heuristic) method to obtain a candidate solution $\bar x$. (iii) Bound the {\em optimality gap} $\phi(\bar x) - \inf \phi$ of $\bar x$ using
\[
\phi(\bar x) - \ninf_{(u,x)\in\reals^m\times\reals^n} \big\{f(u,x) - \langle y,u \rangle  \big\}.
\]
If this quantity is sufficiently close to zero, then $\bar x$ might be deemed acceptable, making further efforts to find an even better solution superfluous. An optimality gap is measured in the units of the objective function and is thus well understood by a decision maker in most cases. We refer to this approach as {\em Rockafellian relaxation}.

Of course, Rockafellian relaxation is only meaningful when the Rockafellian has been selected in such a manner that \eqref{eqn:relax} is easier to solve than the actual problem. This is often the case because  ``lifting'' of the problem to a higher dimension, involving both $x$ and $u$, allows for more flexibility. Moreover, some common choices of Rockafellians lead to explicit expressions for the minimization over $u$ in \eqref{eqn:relax}. Let's denote by
\begin{equation}\label{eqn:lagrangian}
l(x,y) = \ninf_{u\in \reals^m} \big\{f(u,x) - \langle y,u \rangle  \big\}
\end{equation}
the resulting minimum value after such optimization over $u$.  The function $l:\reals^n\times\reals^m\to \Reals$ given by this formula is the {\em Lagrangian}\index{Lagrangian} of $f$. If the Lagrangian corresponding to a particular Rockafellian can be expressed in a convenient form, then the solution of \eqref{eqn:relax} might simplify.

\begin{example}
\label{eLagrCalcul}{\rm (Lagrangian for equalities and inequalities)}. For $f_i:\reals^n\to \reals$, $i=0,1, \dots, m$, and $g_i:\reals^n\to \reals$, $i=1, \dots, q$, let's consider the problem
\[
\nnmin_{x\in\reals^n} \,f_0(x) + \iota_D\big(F(x)\big),
\]
where $D = \{0\}^m\times (-\infty,0]^q$ and
\[
F(x) = \big(f_1(x), \dots, f_m(x), g_1(x), \dots, g_q(x)\big).
\]
A Rockafellian $f:\reals^{m+q}\times\reals^n\to \Reals$ for the problem is defined by
\[
f(u,x) = f_0(x) + \iota_D\big(F(x)+u\big).
\]
The actual problem is then to minimize $f(0,\cdot\,)$ and the Lagrangian has
\[
l(x,y) =\begin{cases}
f_0(x) + \big\langle F(x), y\big\rangle &\mbox{ if } y_{m+1}, \dots, y_{m+q}\geq 0\\
-\infty &\mbox{ otherwise.}
\end{cases}
\]
\end{example}
\state Detail.  With $x\in\reals^n$ and $y\in\reals^{m+q}$, the Rockafellian $f$ produces the Lagrangian given by
\[
l(x,y) = \ninf_{u\in\reals^{m+q}} \big\{ f_0(x) + \iota_D\big(F(x)+u\big) - \langle y, u\rangle\big\}.
\]
If there's  $y_i<0$ for some $i\in \{m+1, \dots, m+q\}$, then we can select $u_j = -f_j(x)$ for all $j\in \{1, \dots, m\}$ and $u_j = -g_j(x)$ for all $j\in \{m+1, \dots, m+q\}\setminus \{i\}$ so that $\iota_D(F(x)+u)$ remains zero as $u_i\to -\infty$. But, then
\[
f_0(x) + \iota_D\big(F(x)+u\big) - \langle y, u\rangle\to -\infty
\]
and $l(x,y) = -\infty$.

If $y_i \geq 0$ for all $i\in \{m+1, \dots, m+q\}$, then $\bar u$, with components $\bar u_j = -f_j(x)$ for all $j\in \{1, \dots, m\}$ and $\bar u_j = -g_j(x)$ for all $j\in \{m+1, \dots, m+q\}$, solves
\[
\nnmin_{u\in \reals^{m+q}} \,f_0(x) + \iota_D\big(F(x)+u\big) - \langle y, u\rangle
\]
and this results in
\[
l(x,y) = f_0(x) - \langle y, \bar u\rangle = f_0(x) + \big\langle y, F(x)\big\rangle.
\]
This establishes the claimed formula for the Lagrangian.\eop

For $h:\reals^m\to \Reals$, we recall that the function $h^*:\reals^m\to \Reals$ defined by
\[
h^*(v) = \nsup_{u\in \reals^m} \big\{ \langle v, u\rangle - h(u)\big\}
\]
is the {\em conjugate} of $h$. This helps us to express Lagrangians for problems in the composite form as seen next; a proof is given by Proposition 5.28 in Royset and Wets \cite{primer}.

\begin{proposition}{\rm (Lagrangian for composite function).}\label{pLagrangianComposite}
For $f_0:\reals^n\to \reals$, $F:\reals^n\to \reals^m$ and proper, lsc and convex  $h:\reals^m\to \Reals$, consider the problem
\[
\nnmin_{x\in X\subset\reals^n} \,f_0(x) + h\big(F(x)\big).
\]
The Rockafellian given by
\[
f(u,x) = \iota_X(x) + f_0(x) + h\big(F(x)+u\big)
\]
recovers the actual problem as minimizing $f(0,\cdot\,)$ and produces a Lagrangian with
\[
l(x,y) = \iota_X(x) + f_0(x) + \big\langle F(x), y\big\rangle - h^*(y).
\]
\end{proposition}

Let's consider a concrete application of Rockafellian relaxation.

\begin{example}\label{eCSPP}{\rm (constrained shortest path problem).}
Let $(V,E)$ be a  directed graph with vertex set $V$ and edge set $E$. Each  edge $(i,j) \in E$ connects distinct  vertices $i, j \in V$,  and it possesses  length $c_{ij} \in [0,\infty)$ and  weights $D_{kij} \in [0,\infty)$ for $k=1, \dots, q$.
A directed $s$-$t$  path is an ordered set of edges of the form $\{ (s,i_1),(i_1,i_2)$,$\dots$, $(i_{\nu-1},t)\}$ for some $\nu\in\nats$. Given two distinct vertices $s,t \in V$, the {\em shortest-path problem} seeks to determine a directed $s$-$t$ path such that the sum of the edge lengths along the path is minimized. This is a well-studied problem that can be solved efficiently using specialized algorithms; see Chapters 4 and 5 in Ahuja et al. \cite{AhujaMagnantiOrlin.93}.

For nonnegative $d_k, k=1, \dots, q$, the task becomes significantly harder if the sum of the weights $D_{kij}$ along the path can't exceed $d_k$ for each $k$. This is the {\em constrained shortest-path problem}, which can be addressed by Rockafellian relaxation---also called Lagrangian relaxation in the  present context. In routing of a drone through a discretized 3-dimensional airspace, the weights $D_{1ij}$ might represent fuel consumption along edge $(i,j)$, which can't exceed a capacity $d_1$. Figure \ref{fig:shortestpath2} illustrates a route satisfying such a fuel constraint while minimizing exposure to enemy radars expressed by $c_{ij}$; cf. Royset et al. \cite{RoysetCarlyleWood.09}.
\end{example}

\drawing{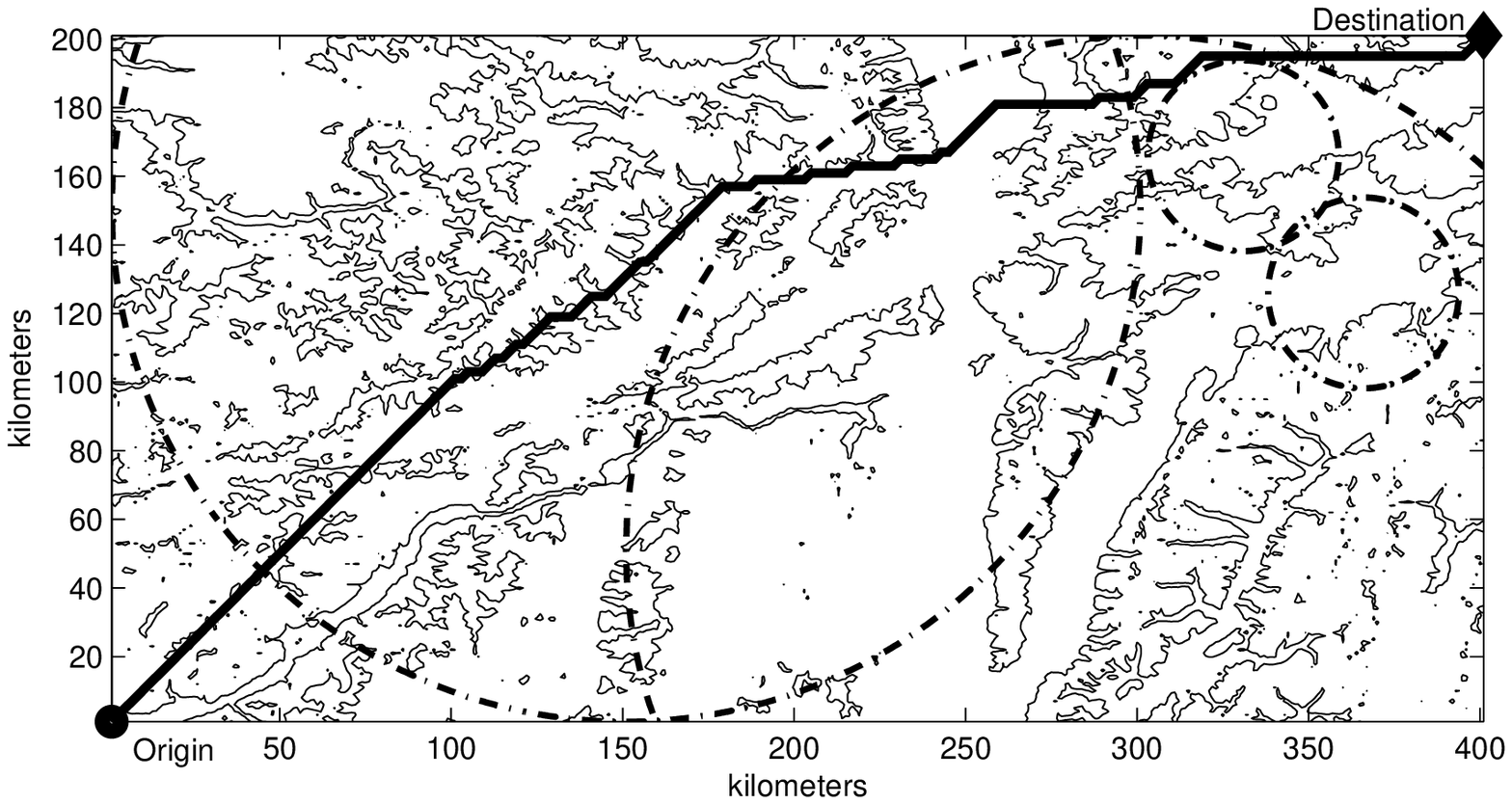}{4.8in}{Route for a drone through a 3-dimensional airspace to a destination (solid line) that minimizes exposure to enemy radars (circles) while satisfying a fuel constraints Royset et al. \cite{RoysetCarlyleWood.09}. Altitude changes to leverage terrain masking aren't shown.}{fig:shortestpath2}

\state Detail.  The constrained shortest-path problem is formulated as follows. Suppose that $m$ is the number of vertices in $V$ and $n$ is the number of edges in $E$. Let $A$ denote the $m \times n$-vertex-edge incidence matrix such that if edge $e = (i,j)$, then $A_{ie} = 1$, $A_{je} = -1$ and $A_{i'e} = 0$ for any $i' \neq i,j$.  Also, let $b_s = 1$, $b_t = -1$ and $b_i = 0$ for $i \in V \backslash \{s,t\}$ and collect them in the vector $b$. For each $k=1, \dots, q$, we place the edge weights $\{D_{kij}, (i,j) \in E\}$ in the vector $D_k$. The $q\times n$-matrix $D$ has $D_k$ as its $k$th row. Let $d=(d_1, \dots, d_q)$. With $c$ being the vector of $\{c_{ij}, (i,j)\in E\}$, the constrained shortest-path problem may then be formulated as (cf. p. 599 of Ahuja et al \cite{AhujaMagnantiOrlin.93})
\[
\nnmin \,\langle c, x\rangle~~\mbox{ subject to } ~A x = b, ~D x \leq d, ~ x \in \{0,1\}^n.
\]
A point $\bar x \in \{0,1\}^n$ satisfying $Ax = b$ corresponds to an $s$-$t$ path with $\bar x_{ij} = 1$ if edge $(i,j)$ is on the path  and $\bar x_{ij} = 0$ otherwise. We assume there's  at least one such path. Then, $\langle c,\bar x\rangle$ gives the length of the path and $\langle D_k, \bar x\rangle$ the $k$th weight of the path.

In the absence of the weight-constraint $D x \leq d$, the model reduces to one for the shortest-path problem and this opens up an opportunity for applying Rockafellian relaxation via Proposition \ref{pLagrangianComposite}. Let
\[
X = \big\{x\in \{0,1\}^n~\big|~Ax = b\big\},
\]
which is nonempty by assumption, and consider the Rockafellian given by
\[
f(u,x) = \iota_X(x) + \langle c, x\rangle +  \iota_{(-\infty,0]^q} (Dx - d+u).
\]
The actual problem corresponds to minimizing $f(0,\cdot\,)$ over $\reals^n$. These definitions fit the setting of Proposition \ref{pLagrangianComposite} and
\[
l(x,y) = \iota_X(x) + \langle c, x\rangle  + \langle Dx - d, y\rangle - \iota_{[0,\infty)^q}(y),
\]
where we use the fact that the conjugate of $\iota_{(-\infty,0]^q}$ is $\iota_{[0,\infty)^q}$; see Example 5.29 of Royset and Wets \cite{primer}. Thus, the chosen Rockafellian results in a Lagrangian without the constraint $Dx\leq d$. The minimization of this Lagrangian with respect to $x$ is nothing but a shortest-path problem on the directed graph, but with edge lengths changed from $c_{ij}$ to $c_{ij} + \nsum_{k=1}^q D_{kij} y_k$, which can be solved efficiently using specialized algorithms. The resulting minimum value, modified by $\langle d, y\rangle$, yields a lower bound on the minimum value of the constrained shortest-path problem as seen from \eqref{eqn:lagrangianLB}. The lower bound can be used to assess the optimality gap for any candidate path, for example obtained by a greedy search or an enumeration algorithm; see Carlyle et al. \cite{CarlyleRoysetWood.08}.\eop

\subsection{Dual Problems}\label{sec:dualproblem}

Rockafellian relaxation produces a lower bound on the minimum value of the actual problem {\em regardless} of the choice of auxiliary vector $y$. However, it would be useful to tune this vector so that the lower bound becomes as high as possible. The underlying model may also be good in the sense that the lower bound associated with the chosen Rockafellian can be brought the whole way up to the minimum value of the actual problem. This could present a major computational advantage: with a properly tuned $y$, minimization of the Rockafellian modified by a linear term can largely substitute for the actual problem. The process of tuning $y$ involves optimization and gives rise to {\em dual problems.}

The notation $y$ for the auxiliary vector in Rockafellian relaxation isn't coincidental. It's deeply connected with multiplier vectors emerging from the Rockafellar condition for optimality; see Theorem \ref{tFermatPara}. In fact, solving a dual problem can be a viable approach to determining a multiplier vector and thus identify sensitivity of the minimum value to perturbations, for example via Proposition \ref{pMultAsSubgrad}.

To make the setting concrete, let $\phi:\reals^n\to \Reals$ and consider the problem
\[
\nnmin_{x\in\reals^n} \,\phi(x),
\]
with an associated Rockafellian $f:\reals^m\times\reals^n\to \Reals$ and anchor at $0$. Via the corresponding Lagrangian $l:\reals^n\times\reals^m\to \Reals$, given by \eqref{eqn:lagrangian}, this produces a {\em dual problem}
\[
\nnmax_{y\in\reals^m} ~\psi(y) = \ninf_{x\in \reals^n} l(x,y).
\]
Since $\inf \phi \geq \psi(y)$ for all $y$ by \eqref{eqn:lagrangianLB}, the dual problem indeed aims to find the best lower bound on $\inf \phi$. Since {\em every} Rockafellian associated with the actual problem defines a dual problem, there are endless possibilities. One might choose a Rockafellian that produces a simple dual problem solvable by standard algorithms. This has been the traditional approach in linear programming. In the setting of Example \ref{eLagrCalcul}, with all functions being affine, the Rockafellian adopted there (focusing on right-hand side perturbations) recovers the usual linear programming dual problems; see Example 5.41 of Royset and Wets \cite{primer} for details. But, this is just one possibility. Beyond the convex case, the dual problems tend to become less tractable but much depends on the structure of the actual problem as well as the choice of Rockafellian. One property is common across all dual problems: the objective function is concave. We realize this by writing
\[
-\psi(y) = -\ninf_{x,u} \big\{f(u,x) - \langle y,u\rangle\big\} = \nsup_{x,u} \big\{\langle y,u\rangle -f(u,x)\big\}.
\]
Thus, $-\psi$ is convex by virtue of being given by the pointwise supremum across a collection of affine functions. Interestingly, regardless of convexity of the actual primal problem and the chosen Rockafellian, the dual objective function $\psi$ is concave, which makes subgradient, cutting plane and proximal point methods at least conceptually available.

Traditionally, the term ``Lagrangian'' has often been limited to the function emerging in the context of equality and inequality constraints under right-hand side perturbations; see Example \ref{eLagrCalcul}.
Even with the broader definition in \eqref{eqn:lagrangian}, we view a Lagrangian as a secondary quantity stemming from a more fundamental Rockafellian. Still, Lagrangians remain important in saddle point theory, which connects the multipliers from optimality conditions with dual variables, serve as a bridge to game theory and promote an elegant symmetry with the actual problem, also expressible in terms of a Lagrangian. We refer to Chapter 5 of Royset and Wets \cite{primer} for details.

In our setting, however, we can completely bypass Lagrangians. Specifically, for any Rockafellian $f:\reals^m \times \reals^n\to \Reals$ and the corresponding Lagrangian $l$, the dual objective function has
\begin{align*}
\psi(y) &= \ninf_{x} l(x, y) = \ninf_{u,x} \big\{f(u,x) - \langle u,y\rangle \big\}\\
& = -\nsup_{u,x} \big\{\langle u,y\rangle + \langle x, 0\rangle - f(u,x)\big\} = -f^*(y,0).
\end{align*}
Consequently, a dual objective function can just as well be defined directly in terms of the conjugate of the underlying Rockafellian.

The usefulness of a dual problem depends on the size of the resulting {\em duality gap}
\[
\inf \phi - \sup \psi.
\]
If it's large, then the corresponding Rockafellian relaxation might not be practically helpful. We say that a Rockafellian possesses {\em strong duality} relative to the actual problem when the duality gap is zero. Then, the actual problem has the good property that its minimum value can be determined by solving a convex optimization problem! We indeed have strong duality in linear programming under the usual right-hand side perturbations provided that the actual problem isn't infeasible. However, strong duality is far from automatic.

\begin{example}{\rm (failure of strong duality).} For the problem of minimizing $x^3$ subject to $x \geq 0$ and the Rockafellian given by
\[
f(u,x) = x^3 + \iota_{(-\infty,0]}(-x+u),
\]
we obtain a Lagrangian of the form
\[
l(x,y) = \begin{cases}
x^3 - xy & \mbox{ if } y\geq 0\\
-\infty & \mbox{ otherwise};
\end{cases}
\]
see Example \ref{eLagrCalcul}. The dual objective function has $\psi(y) = -\infty$ for all $y\in\reals$, while the minimum value of the actual problem is zero so the duality gap is $\infty$.
\end{example}
\state Detail.  In this case, the Lagrangian isn't convex in its first argument. However, strong duality may fail even under convexity. Consider the problem
\[
\nnmin_{x\in\reals^2} \,e^{-x_1} ~\mbox{ subject to }~ g(x) \leq 0, ~\mbox{ where }~ g(x) = \begin{cases}
  x_1^2/ x_2 & \mbox{ if } x_2>0\\
  \infty & \mbox{ otherwise},
\end{cases}
\]
and a Rockafellian of the form
\[
f(u,x) = \begin{cases}
e^{-x_1} & \mbox{ if } g(x) + u \leq 0\\
\infty & \mbox{ otherwise}.
\end{cases}
\]
Similar to Example \ref{eLagrCalcul}, this produces a Lagrangian with
\[
l(x,y) = \begin{cases}
e^{-x_1} + yg(x) & \mbox{ if } x \in \dom g,~ y\geq 0\\
\infty & \mbox{ if } x \not\in \dom g\\
-\infty & \mbox{ otherwise}.
\end{cases}
\]
Consequently, the dual objective function has $\psi(y) = 0$ if $y\geq 0$, but $\psi(y)=-\infty$ otherwise. The maximum value of the dual problem is therefore 0. The actual problem has minimum value of 1. Thus, the duality gap is 1 even though $l(\cdot\,,y)$ is convex regardless of $y\in\reals$.
\eop

Despite these discouraging examples, there are large classes of problems beyond linear optimization problems for which strong duality holds.  The following result summarizes key insights from Theorem 11.39 in Rockafellar and Wets \cite{VaAn} in the convex setting. For strong duality without convexity, we refer to Section 6.B of Royset and Wets \cite{primer}.

Generally, we denote by $\nt C$ the {\em interior} of a set $C\subset\reals^n$, which informally equals $C$ without its boundary points.

\begin{theorem}{\rm (strong duality).}\label{thm:strongduality} For the problem of minimizing $\phi:\reals^n\to \Reals$ and a proper, lsc and convex Rockafellian $f:\reals^m\times\reals^n\to \Reals$ with anchor at $0$, the corresponding dual objective function $\psi$ satisfies strong duality provided that $0 \in \nt(\dom p)$, where $p$ is the minimum value function given by $p(u) = \inf f(u,\cdot\,)$.

If in addition $p(0)>-\infty$, then
\[
\partial p(0) = \nargmax \psi,
\]
which must be a nonempty and bounded set.
\end{theorem}

In Example \ref{eSensitivity}, $\dom p = (-\infty, 1]$ so we certainly have $0 \in \nt(\dom p)$ and strong duality holds; see Figure \ref{fig:sensitivity}. Moreover, $p(0) = 5$ and $\partial p(0) = \{2\}$, which imply that the dual problem has $2$ as its unique maximizer, with 5 as maximum value. We can determine all of this based on Theorem \ref{thm:strongduality} without having a detailed formulation for the dual problem. The formula for $\partial p(0)$ supplements Theorem \ref{tParametric}, but most significantly it highlights the profound role played by a dual problem. Under the conditions of the theorem, solving the dual problem furnishes both the minimum value of the actual problem as well as its sensitivity to perturbations as defined by a Rockafellian.

In a specific setting, strong duality is guaranteed by ensuring that the constraint functions leave some ``slack,'' which often is easily verified.

\begin{example}{\rm (Slater constraint qualification).}\label{eSlater} For smooth convex functions $f_0,g_i:\reals^n\to \reals$, $i=1, \dots, q$ and the problem
\[
\nnmin_{x\in\reals^n} \,f_0(x) ~\mbox{ subject to }~ g_i(x) \leq 0, ~i=1, \dots, q,
\]
let's consider the Rockafellian given by
\[
f(u,x) = f_0(x) + \iota_{(-\infty, 0]^q} \big(G(x) + u \big), ~\mbox{ with } G(x) = \big(g_1(x), \dots, g_q(x)\big).
\]
The resulting dual problem satisfies strong duality provided that the following  {\em Slater condition} holds:
\[
\exists \bar x ~\mbox{ such that } ~g_i(\bar x) < 0, ~~i=1, \dots, q.
\]
\end{example}
\state Detail.  Let $p(u) = \inf f(u,\cdot\,)$. Under the Slater condition, there exist
$\bar x\in\reals^n$ and $\delta>0$ such that $g_i(\bar x) + u_i \leq 0$ when $|u_i| \leq \delta$ for all $i$. Consequently, $p(u) \leq f_0(\bar x)\in\reals$ when $\|u\|_\infty \leq \delta$, which means that $0\in \nt (\dom p)$ and strong duality holds by Theorem \ref{thm:strongduality}.

In the setting of Theorem 7.2, specialized to convex inequality constraints, the Slater condition ensures that the qualification \eqref{eqn:compositeQual} holds; see Examples 4.49 and 5.47 of Royset and Wets \cite{primer} for details.

The Slater condition is by no means necessary for strong duality. For example, consider the problem of minimizing $x$ subject to $x^2 \leq 0$ and the Rockafellian given by
\[
f(u,x) = x + \iota_{(-\infty, 0]}(x^2+u).
\]
We obtain directly that
\[
p(u) = \inf f(u,\cdot\,) = \begin{cases}
-\sqrt{-u} & \mbox{ for } u\leq 0\\
\infty & \mbox{ otherwise}
\end{cases}
\]
so the requirement $0\in \nt(\dom p)$ of the Strong Duality Theorem \ref{thm:strongduality} doesn't hold.
The corresponding Lagrangian has
\[
l(x,y) = \begin{cases}
x + yx^2 & \mbox{ for } y\geq 0\\
-\infty & \mbox{ otherwise}
\end{cases}
\]
by Example \ref{eLagrCalcul} and the dual objective function has
\[
\psi(y) = \begin{cases}
-1/(4y) & \mbox{ for } y>0\\
-\infty & \mbox{ otherwise}.
\end{cases}
\]
Consequently, $p(0) = \nsup \psi = 0$ so strong duality holds, but the Slater condition fails.\eop

The requirement $0 \in \nt(\dom p)$ in the Strong Duality Theorem \ref{thm:strongduality} is ensured when $f(u, \cdot\,)$ approximates $f(0, \cdot\,)$ in a certain sense for $u$ near zero. This is formalized in the next statements, which highlight the role of suitable approximations to ensure strong duality; see Theorem 5.49 and Corollary 5.50 in Royset and Wets \cite{primer}.

\begin{theorem}{\rm (strong duality from epi-convergence).}\label{tStrongDualFromEpi}
For the problem of minimizing $\phi:\reals^n\to \Reals$ and a proper Rockafellian $f:\reals^m\times \reals^n\to \Reals$ with anchor at $0$, suppose that there are $u^\nu\to 0$ and $y^\nu\in\reals^m$ such that the following hold:
  \begin{enumerate}[{\rm (a)}]
  \item $f(u^\nu,\cdot\,)\eto f(0,\cdot\,)$
  \item there is a compact set $B$ such that $B \cap \nargmin f(u^\nu,\cdot\,)$ is nonempty for all $\nu$.
  \item $\nliminf\, \langle y^\nu,u^\nu\rangle \leq 0$
  \item $\inf f(u^\nu, \cdot\,) = \sup \psi^\nu = \psi^\nu(y^\nu)$,
  \end{enumerate}
where
\[
\psi^\nu(y) = \ninf_x l^\nu(x,y)~ \mbox{ and }~ l^\nu(x,y) = \ninf_u \big\{ f(u^\nu + u,x) - \langle y, u\rangle\big\}.
\]
Let $\psi$ be the dual objective function produced by $f$ via \eqref{eqn:lagrangian}. Then,
\[
\inf \phi = \sup \psi>-\infty.
\]
\end{theorem}

The theorem shows that if we can construct perturbed functions via a Rockafellian and they epi-convergence to the actual objective function as well as possess a strong duality property, then the resulting dual problem indeed reproduces the minimum value of the actual problem provided that assumption (c) also holds.

Since $u^\nu\to 0$, assumption (c) certainly holds when $\{y^\nu,\nu\in\nats\}$ is bounded. One can view $\psi^\nu$ as a dual objective function produced by the Rockafellian of the form $f^\nu(u,x) = f(u^\nu+u,x)$. The vector $y^\nu$ then solves the dual problem associated with $f^\nu$. Thus, it's plausible that $\{y^\nu,\nu\in\nats\}$ could be bounded.

In some cases, assumption (c) is automatic even when $\{y^\nu,\nu\in\nats\}$ is unbounded. For example, if
\[
f(u,x) = f_0(x) + \iota_{(-\infty, 0]^m} \big(F(x) + u\big)
\]
for $f_0:\reals^n\to \reals$ and $F:\reals^n\to \reals^m$, both smooth, then
\[
\psi^\nu(y) = \ninf_{x\in\reals^n} f_0(x) + \big\langle F(x)+u^\nu, y\big\rangle - \iota_{[0,\infty)^m}(y)
\]
by Proposition \ref{pLagrangianComposite}. Thus, a maximizer $y^\nu$ of $\psi^\nu$ is necessarily nonnegative. We can then choose $u^\nu \leq 0$ so that $\langle y^\nu,u^\nu\rangle \leq 0$ and  $f(u^\nu,\cdot\,)\eto f(0,\cdot\,)$, which can be seen by working directly from Proposition \ref{tEpiCnvr}.

The theorem isn't restricted to any particular type of Rockafellian.  Still, in the convex case, several aspects simplify.

\begin{corollary}
For the problem of minimizing $\phi:\reals^n\to \Reals$ and a proper, lsc and convex Rockafellian $f:\reals^m\times \reals^n\to \Reals$ with anchor at $0$, suppose that there's a compact set $B\subset\reals^n$ and $u^\nu \in \nt(\dom p)\to 0$ such that $f(u^\nu,\cdot\,)\eto f(0,\cdot\,)$ and $B\cap\nargmin f(u^\nu,\cdot\,)$ is nonempty for all $\nu$, where $p$ is the minimum value function given by $p(u) = \inf f(u, \cdot\,)$. Let $\psi$ be the dual objective function produced by $f$ via \eqref{eqn:lagrangian}.

If \,$\inf \phi < \infty$, then
\[
\inf \phi = \sup \psi
\]
and this value is finite.
\end{corollary}

\begin{example}{\rm (Slater constraint qualification; cont.).}\label{eSlater3} The corollary confirms the strong duality assertion towards the end of Example \ref{eSlater} even though the Slater condition fails.
\end{example}
\state Detail.  In this case, the Rockafellian, given by $f(u,x) = x + \iota_{(-\infty, 0]}(x^2+u)$, is proper, lsc and convex. Moreover, $p(u) = -\sqrt{-u}$ for $u\leq 0$ and $p(u) = \infty$ otherwise. Thus, one can take $u^\nu = -1/\nu$ in the corollary and then $f(u^\nu,\cdot\,)\eto f(0,\cdot\,)$; the compactness condition holds since $\dom f(u^\nu, \cdot\,)\subset \dom f(u^1, \cdot\,)$.
\eop\\

\state Acknowledgements. This work is supported in part by the Office of Naval Research under MIPR N0001421WX01496 and the Air Force Office of Scientific Research under MIPR F4FGA00350G004.\\

\bibliographystyle{plain}
\bibliography{refs}

\end{document}